\documentclass{article}
\usepackage{latexsym}
\usepackage{amsfonts}
\usepackage{algorithm}
\usepackage{algpseudocode}
\usepackage{color}
\usepackage{latexsym}
\usepackage{mathtools}
\usepackage{pgf,tikz}
\usetikzlibrary{arrows}
\usepackage{fancybox}
\usepackage{float}
\usetikzlibrary{shapes,arrows}
\usepackage{psfrag}
\usepackage{hyperref}
\newtheorem{proposition}{Proposition}[section]
\newtheorem{assumption}[proposition]{Assumption}
\newtheorem{theorem}[proposition]{Theorem}
\newtheorem{definition}[proposition]{Definition}
\newtheorem{corollary}[proposition]{Corollary}
\newtheorem{lemma}[proposition]{Lemma}
\newtheorem{remark}{Remark}
\newcommand{\proof}{{\bf Proof}}

\newcommand{\commento}[1]{}

\def\eps{\epsilon}

\numberwithin{equation}{section}

\newcommand\Algphase[1]{%
	\Statex\hspace*{-\algorithmicindent}\textbf{#1}%
}

\begin{document}
\pagestyle{empty}

\vskip 2cm \begin{center}{\huge Derivative-free methods for mixed-integer nonsmooth constrained optimization}\end{center}
\par\bigskip
\centerline{\Large 
	T. Giovannelli\footnote{Lehigh University, Department of Industrial and Systems Engineering, (email) {\tt tog220@lehigh.edu}}, 
	G. Liuzzi\footnote{``Sapienza" Universit\`a  di Roma, Dipartimento di Ingegneria Informatica Automatica e Gestionale ``A. Ruberti'', (email) {\tt liuzzi@diag.uniroma1.it}, {\tt lucidi@diag.uniroma1.it}}, 
	S. Lucidi$^{2}$, 
	F. Rinaldi\footnote{Universit\`a  di Padova, Dipartimento di Matematica ``Tullio Levi-Civita'', (email) {\tt rinaldi@math.unipd.it}}
}
\par\bigskip\bigskip

\bigskip

{\small \centerline{\bf Abstract}
\begin{quote}
In this paper, we consider mixed-integer nonsmooth constrained optimization problems whose objective/constraint functions are available
only as the output of a black-box zeroth-order oracle (i.e., an oracle that does not
provide derivative information) and  we propose a new derivative-free linesearch-based algorithmic framework to suitably handle those problems. We first  describe a scheme for bound constrained problems that combines
a dense sequence of directions
(to handle the nonsmoothness of the objective function) with 
primitive directions (to handle discrete variables). Then, we embed an exact penalty approach in the scheme to suitably
manage nonlinear  (possibly nonsmooth) constraints. We analyze the global convergence properties of the proposed algorithms toward stationary points and we report the results of an extensive numerical experience on a set of mixed-integer test problems.
\end{quote} }

\par\medskip


\newpage
\pagestyle{plain}


\section{Introduction}

\noindent We consider the following mixed-integer nonlinearly constrained problem
\par\medskip \noindent
\begin{equation}\label{probconstr_1}
\begin{array}{l}
\min \,\,\, f(x) \\ [.3em] s.t. \quad g(x) \le 0,  \\ [.3em] \,\,\qquad l\le x\le u, \\ [.3em] \qquad \ x_i\in \mathbb{R} \ \text{ for all } \ i\in I^c,
\\ [.3em] \qquad \ x_i\in \mathbb{Z} \ \text{ for all } \ i\in I^z,\end{array}
\end{equation}
\par\medskip \noindent where $x\in  \mathbb{R}^n$, \, $l,u\in  \mathbb{R}^n$, \, and $I^c\cup I^z = \{1,\ldots ,n \}$, with $I^c \cap I^z = \emptyset$. 
We assume $l_i < u_i$ for all $i \in I^c\cup I^z$, and
$l_i,u_i\in  \mathbb{Z}$ for all $i\in I^z$. Moreover, the functions $f: \mathbb{R}^n\to  \mathbb{R}$ and $g: \mathbb{R}^n\to  \mathbb{R}^m$, which may be nondifferentiable, are supposed to be Lipschitz continuous with respect to $x_i$ for all $i\in I^c$, i.e., for all $x,y\in\mathbb{R}^n$ a constant $L>0$ exists such that
\begin{equation}\label{def:Lipschitz}
|f(x)-f(y)| \leq L \|x-y\|,\quad\mbox{with}\ x_i=y_i, \mbox{ for all } \ i\in I^z.
\end{equation}
We define the sets
\begin{eqnarray*}
	&& X:=\{x \in \mathbb{R}^n: l\le x\le u\},\quad {\cal F} := \{x\in \mathbb{R}^n: g(x) \leq 0\},\\
	&& {\cal Z}:=\{x\in \mathbb{R}^n: x_i\in \mathbb{Z} \text{ with } i\in I^z \},
\end{eqnarray*}
and assume throughout the paper that $X$ is a compact set.  Therefore, $l_i$ and $u_i$ cannot be infinite. Problem \eqref{probconstr_1} can hence be reformulated as follows
\begin{equation}\label{probconstr}
\begin{array}{l}
\min f(x) \\ [.3em] s.t. \quad x \in {\cal F} \cap {\cal Z} \cap X.\end{array}
\end{equation}
\par
The objective and constraint functions in \eqref{probconstr} are assumed to be of black-box zeroth-order type (i.e., the analytical expression is unknown, and the function value corresponding to a given point is the only available information). We thus consider black-box Mixed-Integer Nonlinear Programs (MINLPs), a class of challenging problems frequently arising in real-world applications. Those problems are usually solved through tailored derivative-free optimization algorithms (see, e.g.,  \cite{audet.2017,boukouvala.2015,conn.2009,larson.2019} and references therein for further details) able to properly manage the presence of both continuous and discrete variables. 
\par
The optimization methods for black-box MINLPs that we consider in here can be divided into two main classes: direct search and model-based methods. 
The direct search methods for MINLPs usually share two main features: they perform an alternate minimization between continuous and discrete variables, and use a fixed neighborhood to explore the integer lattice. 
In particular, \cite{audet.2001} adapts the Generalized Pattern Search (GPS), proposed in \cite{torczon.1997}, to solve problems with categorical variables (those variables include integer variables as a special case), so-called mixed-variable problems. 
The approach in \cite{audet.2001} has been then extended to address problems with general constraints~\cite{abramson.2004} and stochastic objective function~\cite{sriver.2009}. In \cite{abramson.2004}, constraints are tackled by using a  filter approach similar to the one described in~\cite{audet.2004}. 
Derivative-free methods for categorical variables and general constraints  have been also studied in \cite{lucidi.2005} and \cite{abramson.2009}. In particular, \cite{lucidi.2005} proposes a general algorithmic framework whose global convergence holds for any continuous local search (e.g., a pattern search) satisfying suitable properties. In \cite{abramson.2009}, the class of Mesh Adaptive Direct Search (MADS), originally introduced in \cite{audet.2006} to face nonsmooth nonlinearly constrained problems, is extended to solve mixed-variable problems. 
 Constraints are tackled through an extreme barrier approach in this case. The original MADS algorithm has been recently extended in \cite{audet.2019} to solve problems with granular variables, i.e., variables with a controlled number of decimals (integer variables can be seen as a special case as well), and the objective function is nonsmooth over the continuous variables. 
In addition to these references, another work that is worth mentioning is \cite{porcelli.2017}, where a mesh-based direct search algorithm is proposed for bound constrained mixed-integer problems whose objective function is allowed to be nonsmooth and noncontinuous.
\par
In \cite{liuzzi.2012}, three algorithms are proposed for bound constrained MINLP problems. Differently from the aforementioned works, the discrete neighborhood does not have a  fixed structure, but depends on a linesearch-type procedure. 
The first algorithm in \cite{liuzzi.2012} was extended in \cite{liuzzi.2014} and \cite{yang.2019}: the former deals with the constrained case by adopting a sequential penalty approach, while the latter replaces the maximal positive basis with a minimal positive basis based on a direction-rotation technique. Bound constrained MINLP problems are considered also in \cite{garcia-palomares.2012}, which extends the algorithm for continuous smooth and nonsmooth objective functions introduced in \cite{garcia-palomares.2002}.  
\par
We report here some other direct search methods (not directly connected with MINLP problems) that is worth  mentioning for their influence on algorithm development. In \cite{fasano.2014}, the authors propose a new linesearch-based method for nonsmooth nonlinearly constrained optimization problems, ensuring convergence towards Clarke-Jahn stationary points. The constraints are tackled through an exact penalty approach. 
In \cite{custodio.2007} and \cite{custodio.2008}, the authors analyze the benefit in terms of efficiency deriving from different ways of incorporating the simplex gradient into direct search algorithms (e.g., GPS and MADS) for minimizing objective functions which not necessarily require continuous differentiability. 
In \cite{vicente.2012}, the authors analyze the convergence properties of direct search methods applied to the minimization of discontinuous functions. 
\par
Model-based methods are also widely used in this context. In \cite{costa.2018}, the authors describe an open-source library, called RBFOpt,  that uses surrogate models based on radial basis functions
for handling bound constrained MINLPs.
The same class of problems is  also tackled in \cite{newby.2014} through quadratic models. This paper extends to the mixed-integer case the trust-region derivative-free algorithm BOBYQUA introduced in \cite{powell.2009} for continuous problems.   
Surrogate models employing radial basis functions are used in \cite{muller.2013a} to describe an algorithm, called SO-MI, able to converge to the global optimum almost surely. 
A similar algorithm, called SO-I, is proposed by the same authors in \cite{muller.2013b} to address integer global optimization problems. 
In \cite{muller.2015}, the authors propose an algorithm for MINLP problems that modifies the sampling strategy used in SO-MI and uses also an additional local search. Finally, Kriging models were effectively used in \cite{hemker.2008} and \cite{halstrup.2016} to develop new sequential algorithms.
Models can also be used to boost direct-search methods. For example, in NOMAD, i.e., the software package that implements the MADS algorithm (see \cite{Nomad}), a surrogate-based model is used to generate promising points.

\par
In \cite{larson.leyffer.2019} and \cite{liuzzi.2020} methods for black-box problems with only unrelaxable integer variables
are devised. In particular, the authors in \cite{larson.leyffer.2019} propose a method for minimizing convex black-box integer problems that uses secant functions interpolating previous evaluated points. In \cite{liuzzi.2020}, a new method based on a nonmonotone linesearch and primitive directions is proposed to solve a more general problem where the objective function is allowed to be nonconvex. The primitive directions allow the algorithm to escape bad local minima, thus providing the potential to find a global optimum (this anyway requires the exploration of large neighborhoods).

In this paper, we propose new derivative-free linesearch-type algorithms for mixed-integer nonlinearly constrained problems with possibly nonsmooth functions. We combine the strategies successfully tested in \cite{fasano.2014} and \cite{liuzzi.2020} for continuous and integer problems, respectively, to
devise a globally convergent 
algorithmic framework that
enables us to tackle the mixed-integer case.
Continuous and integer variables are suitably handled by means of specific local searches in this case. On the one hand, a dense sequence of search directions is used to explore the subspace related to the continuous variables 
and detect descent directions, whose cone can be arbitrarily narrow due to nonsmoothness. On the other hand, a set of primitive discrete directions is adopted to guarantee a thorough exploration of the integer lattice in order to escape bad local minima.  We develop a first algorithm for bound constrained problems, then we adapt it to handle the presence of general nonlinear constraints by using an exact penalty approach. Since only the violation of such constraints is included in the penalty function, the algorithm developed for bound constrained problems can be easily adapted to minimize the penalized problem. 

With regard to the convergence results, we can prove that particular sequences of iterates yielded by the two algorithms converge to suitably defined stationary points of the problem considered. In the generally constrained case, this result is based on the equivalence between the original problem and the penalized problem. 

The paper is organized as follows. In Section~\ref{sec:def}, we report some definitions and preliminary results. In Section~\ref{sec:bound_constrained_case}, we describe the algorithm proposed for mixed-integer problems with bound constraints and we analyze its convergence properties. The same type of analysis is reported in Section~\ref{sec:nonlinearly_constrained_case} for the algorithm addressing mixed-integer problems with general nonlinear constraints. Section~\ref{sec:numerical_exp} describes the results of extensive numerical experiments performed for both algorithms. Finally, in Section~\ref{sec:conclusion} we include some concluding remarks and we discuss future work.

\section{Notation and preliminary results}\label{sec:def}
Given a vector $v\in  \mathbb{R}^n$, we introduce the subvectors $v_c\in  \mathbb{R}^{|I^c|}$ and $v_z\in
	\mathbb{R}^{|I^z|}$, given by
	$$v_c=\left[v_i\right]_{i\in I^c}\quad \text{and}\quad v_z=\left[v_i\right]_{i\in I^z},$$
	where $v_i$ denotes the i-th component of $v$. When a vector is an element of an infinite sequence of vectors $\{v_k\}$, the i-th component will be denoted as $(v_k)_i$, in order to avoid possible ambiguities. Moreover, throughout the chapter we denote by $\|\cdot\|$ the Euclidean norm.
\par
The search directions considered in the algorithms proposed in the next sections have either a null continuous subvector or a null discrete subvector, meaning that we do not consider directions that update both continuous and discrete variables simultaneously. We first report the definition of primitive vector, used to characterize the subvectors of the search directions related to the integer variables. Then we move on to the properties of the subvectors related to the continuous variables.

{}From \cite{liuzzi.2020} we report the following definition of primitive vector.

\begin{definition}[Primitive vector]
	\label{def:primitive_vector}
	A vector $v\in\mathbb{Z}^n$ is called primitive if the greatest common divisor of its components $\{v_1,\dots,v_n\}$ is equal to 1.
\end{definition}
\par 
Since the objective and constraint functions of the problem considered are assumed to be nonsmooth (when fixing the discrete variables), proving convergence to a stationary point requires particular subsequences of the continuous subvectors of the search directions to be provided with the density property.  Since the feasible descent directions can form an arbitrarily narrow cone (see, e.g.,  \cite{abramson.2009.orthomads} and \cite{audet.2006}), a finite number of search directions is indeed not sufficient. Denoted the unit sphere with center in the origin as $$S(0,1)=\{s \in \mathbb{R}^n \ : \ \|s_c\| = 1 \text{ and } \|s_z\| = 0\},$$ we extend to the mixed-integer case the definition of a dense subsequence of directions.
\begin{definition}[Dense subsequence]\label{def:dense_sequence}
	Let $K$ be an infinite subset of indices (possibly $K=\{0,1,\dots\}$) and $\{s_k\} \subset S(0,1)$ a given sequence of directions. The subsequence $\{s_k\}_K$ is said to be dense in  $S(0,1)$ 
	if, for any  $\bar s\in S(0,1)$ and for any $\eps > 0$, there exists an index $k\in K$ such that $\|s_k-\bar s\|\leq \eps$.
\end{definition}

Similarly to what is done in \cite{abramson.2009}, we extend to the mixed-integer case the definition of generalized directional derivative, which is also called Clarke directional derivative, given in~\cite{clarke.1983}. This allows us to provide necessary optimality conditions for Problem~\eqref{probconstr}. We also recall the definition of generalized gradient.
\begin{definition}[Generalized directional derivative and generalized gradient] \label{def:clarke_directional_derivative}
	Let $h: \mathbb{R}^n \to \mathbb{R}$ be a Lipschitz continuous function near $x \in \mathbb{R}^n$ with respect to its continuous variables $x_c$ (see, e.g., (\ref{def:Lipschitz})). The generalized directional	derivative of $h$ at $x$ in the direction $s \in \mathbb{R}^n$, with $s_i = 0$ for $i \in I^z$, is 
	\begin{equation}\label{eq:clarke_directional_derivative}
	h^{Cl}_{x_c}(x; s) = \limsup_{\footnotesize\begin{array}{l}y_c\to x_c, y_z = x_z, t\downarrow 0\end{array}} \frac{h(y+t s) -
		h(y)}{t}.
	\end{equation}
	To simplify the notation, the generalized gradient of $h$ at $x$ w.r.t the continuous variables can be redefined as
	$$\partial_{x_c} h( x)=\{v\in \mathbb{R}^{n} :\ v_i=0,\ i\in I^z,\ \mbox{and}\  h_{x_c}^{Cl}( x; s)\ge s^Tv \ \text{ for all } s \in \mathbb{R}^n, \text{ with } s_i = 0 \text{ for } i \in I^z \}.$$  
\end{definition}
\smallskip
\par\noindent
Moreover, let us denote the orthogonal projection over the set $X$ as $[x]_{[l,u]}=\max\{l,\min\{u,x\}\}$ and the interior of a set $\cal C$ as $\stackrel{\circ}{\cal C}$. These concepts will be used throughout the paper.

\subsection{The bound constrained case}\label{subsec:def_bound_constr}
We first focus on a simplified version of Problem~\eqref{probconstr_1}, where we only
consider bound constraints in the definition of the feasible set:
\par\medskip \noindent
\begin{equation*}\label{probbox_1}
\begin{array}{l}
\min \,\,\, f(x) \\ [.3em] s.t. \quad l\le x\le u, \\ [.3em] \qquad \ x_i\in \mathbb{R} \ \text{ for all } \ i\in I^c,
\\ [.3em] \qquad \ x_i\in \mathbb{Z} \ \text{ for all } \ i\in I^z.\end{array}
\end{equation*}
\par\medskip \noindent Such a problem can be reformulated as follows
\begin{equation}\label{probbox}
\begin{array}{l}
\min f(x) \\ [.3em] s.t. \quad x \in X \cap {\cal Z}.\end{array}
\end{equation}
\par\medskip \noindent
The next definitions are related to directions that are feasible with respect to $X \cap {\cal Z}$.
\begin{definition}[Set of feasible primitive  directions]\label{def:set_feasible_primitive_discrete_directions}
	Given a point $ x\in X\cap {\cal Z}$, 
	\begin{eqnarray*}
		D^z( x) & = &\{d\in\mathbb{Z}^n: d_z \ \hbox{is a primitive vector},\ d_i=0 \hbox{ for all } i\in I^c, \, and\\
		& & \qquad\qquad\	x+ d \in X\cap\cal{Z}\}
	\end{eqnarray*}
is the set of feasible primitive  directions at $x$ with respect to $X\cap\cal{Z}$.
\end{definition}
\begin{definition}[Union of the sets of feasible primitive directions]\label{def:union_set_feasible_primitive_discrete_directions} 
	\begin{eqnarray*}
		\bar D & = & \bigcup_{x \in X\cap\cal{Z}} D^z(x)
	\end{eqnarray*}
	is the union of the sets of feasible primitive  directions with respect to $X\cap\cal{Z}$.
\end{definition}
\begin{proposition}\label{prop:union_set_feasible_discr_dir} 
	The union of the sets of feasible primitive discrete directions $\bar D$ has a finite number of elements.
\end{proposition}
\proof \quad
Given $d \in \bar D$, it follows that $d \in D^z(x)$ for some $x \in X \cap \mathcal{Z}$. By the definition of $D^z(x)$, we have that $d_i=0$ for all $i \in I_c$ and $d_z$ is a primitive vector. Hence, by considering the boundedness of $X$, for all $x \in X\cap\mathcal{Z}$ the number of subvectors $d_z, \ d \in D^z(x)$, is finite. By the boundedness of $X$ and by $d_i=0 \in I_c$, it follows that the number of distinct $x_z$ is finite and that the number of distinct $D^z(x),$ with $x \in X\cap\mathcal{Z}$, is finite as well. Therefore, the union of subvectors $d_z$ from a finite number of distinct $D^z(x),$ with $x \in X \cap \mathcal{Z}$, has a finite number of elements. 
$\hfill\Box$ \\
\par
The cone of feasible continuous directions is defined according to the following definition.
\begin{definition}[Cone of  feasible continuous directions]\label{def:cone_feasible_continuous_directions} Given a point $x\in X\cap\cal{Z}$, the set
	\begin{eqnarray*}
		D^c(x) & = &\{s\in \mathbb{R}^n: s_i=0\quad \hbox{for all}\quad i\in I^z,\\
		&   &\hspace*{1.5cm} s_i\ge 0 \quad \hbox{for all}\quad  i\in I^c\quad\mbox{and}\quad x_i = l_i,\\
		& &\hspace*{1.5cm}  s_i\le 0\quad \hbox{for all}\quad  i\in I^c\quad\mbox{and}\quad  x_i = u_i, \\ 
		&   & \hspace*{1.5cm} s_i\in \mathbb{R}\quad \hbox{for all}\quad  i\in I^c\quad\mbox{and}\quad l_i<x_i < u_i,\}
	\end{eqnarray*}
	is the cone of feasible continuous directions at $x$ with respect to $X\cap\cal{Z}$.
\end{definition}
\par
Now we report a technical proposition whose proof can be easily derived from the proof of \cite[Proposition 2.3]{lin.2009}.
\begin{proposition}\label{deramm}
	Let $\{x_k\} \subset X\cap\cal{Z}$ for all $k$, and $\{x_k\}\to\bar x\in  X\cap\cal{Z}$ for $k\to\infty$. Then, for $k$
	sufficiently large,
	$$D^c(\bar x)\subseteq D^c(x_k).$$
\end{proposition}
Moreover, we introduce two definitions of neighborhood related to the discrete variables and we recall the definition of neighborhood related to the continuous variables. 
\begin{definition}[Discrete neighborhood]\label{def:discrete_neighborhood}
	Given a point $\bar x\in X\cap\cal{Z}$, the discrete neighborhood of $\bar x$ is
	\[
	{\cal B}^z(\bar x) = \{x\in X\cap{\cal Z}:\quad x = \bar x+d,\quad  \mbox{with}\quad d\in D^z(\bar x)\}.
	\]
\end{definition}
\begin{definition}[Weak discrete neighborhood]\label{def:weak_discrete_neighborhood}
	Given a point $\bar x\in X\cap{\cal Z}$ and a subset $\tilde D^z\subset D^z(\bar x)$,  the weak discrete neighborhood of $\bar x$ is
	\[
	{\cal B}^z_w(\bar x, \tilde D^z) = \{x\in X\cap{\cal Z}:\quad x = \bar x+ d,\quad \mbox{with}\quad d\in \tilde D^z\}.
	\]
\end{definition}
\begin{definition}[Continuous neighborhood]\label{def:continuous_neighborhood}
	Given a point $\bar x\in X\cap\cal{Z}$  and a scalar $\rho$, the continuous neighborhood of $\bar x$ is
\begin{eqnarray*}
 {\cal B}^c(\bar x;\rho)   & = & \left\{ x\in  X: x_z = \bar x_z \ and \ \|x_c-\bar x_c\| \leq \rho\right\}.
\end{eqnarray*}
\end{definition}
%
%

Two definitions of local minimum points are given below.
\begin{definition}[Local minimum point]\label{def:local_minimum_point}
A point $x^*\in X\cap {\cal Z}$ is a local minimum point of Problem (\ref{probbox}) if, for some $\eps >0$,
\begin{eqnarray}
 && f(x^*)\leq f(x) \quad \text{ for all } x\in {\cal B}^c(x^*;\eps), \label{eq:local_minimum_point_cont}\\[1.0em]
 &&f(x^*)\leq f(x) \quad \text{ for all } x\in {\cal B}^z(x^*) \label{eq:local_minimum_point_disc}.
\end{eqnarray}
\end{definition}
\begin{definition}[Weak local minimum point]\label{def:weak_local_minimum_point}
	A point $x^*\in X\cap {\cal Z}$ is a weak local minimum point of Problem (\ref{probbox}) if, for some $\eps >0$ and some $\tilde D^z\subset D^z(x^*)$,
	\begin{equation}\label{eq:weak_local_minimum_point}
	\begin{array}{ll}
	f(x^*)\leq f(x) & \text{ for all } x\in {\cal B}^c(x^*;\eps), \\[1.0em]
	f(x^*)\leq f(x) & \text{ for all } x\in {\cal B}^z_w(x^*, \tilde D^z).
	\end{array}
	\end{equation}
\end{definition} 

We now extend to the mixed-integer case the definition of Clarke-Jahn generalized directional derivative given in \cite[Section 3.5]{jahn.2014}. As opposed to the Clarke directional derivative, in this definition the limit superior is considered only for points $y$ and $y + ts$ in $X\cap{\cal Z}$, thus requiring stronger assumptions.
\begin{definition}[Clarke-Jahn generalized directional derivative] \label{def:clarke-jahn_directional_derivative}
	Given a point $x \in X\cap{\cal Z}$ with continuous subvector $x_c$, the Clarke-Jahn generalized	directional derivative of function $f$ along direction $s \in D^c(x)$ is given by:
	\begin{equation}\label{eq:clarke-jahn_directional_derivative}
	f^\circ_c(x;s) = \limsup_{\footnotesize\begin{array}{l}y_c\to x_c,y_z=x_z, y\in X\cap{\cal Z}\\ \quad \ t\downarrow 0, y+ts\in X\cap{\cal Z}\end{array}}
	\frac{f(y+t s) -f(y)}{t}
	\end{equation} 
\end{definition}

We finally report a few basic stationarity definitions.


	\begin{definition}[Stationary point]\label{def:stationary_point}
		A point $x^*\in X\cap{\cal Z}$ is a stationary point of Problem (\ref{probbox}) when 
\begin{align}
&f^{\circ}_c(x^*; s)\geq 0,\quad\quad \mbox{for all}\ s \in D^c(x^*),\label{nulgrad_unc} \\ 
&f(x^*)\le f(x),\quad\quad \mbox{for all}\ x\in {\cal B}^z(x^*). \label{min_discr_unc}
\end{align}
	\end{definition}
	\begin{definition}[Weak stationary point]\label{def:weak_stationary_point}
		A point $x^*\in X\cap {\cal Z}$ is a weak stationary point of Problem (\ref{probbox}) when, for some $\tilde D^z\subset D^z(x^*)$, 
	\begin{align}
	&f^{\circ}_c(x^*; s)\geq 0,\quad\quad \mbox{for all}\ s \in D^c(x^*), \label{wnulgrad_unc} \\  
	&f(x^*)\le f(x),\quad\quad \mbox{for all}\  x\in {\cal B}^z_w(x^*;\tilde D^z). \label{wmin_discr_unc}
	\end{align}

	\end{definition}
	\begin{definition}[Clarke stationary point]\label{def:clarke_stationary_point}

A point $x^*\in X\cap{\cal Z}$ is a Clarke stationary point of Problem (\ref{probbox}) when it satisfies
		\begin{align}
			&f^{Cl}_c(x^*; s)\geq 0\quad\quad \mbox{for all}\ s \in D^c(x^*),\label{nulgrad_unc_clarke} \\ 
			&f(x^*)\le f(x)\quad\quad \mbox{for all}\ x\in {\cal B}^z(x^*). \label{min_discr_unc_clarke}
		\end{align}
	\end{definition}


	\begin{corollary} \label{cor:clarke_stationary_point}
	Any (weak) minimum point of Problem (\ref{probbox}) is a (weak) stationary point. Furthermore, 
		any stationary point for Problem (\ref{probbox}) is a Clarke stationary point.		
	\end{corollary}
	
\subsection{The nonsmooth nonlinearly constrained case}\label{subsec:def_general_constr}
Now we turn our attention to the general case defined through Problem~\eqref{probconstr}. A local minimum point for this problem is defined as follows.
\begin{definition}[Local minimum point]\label{def:local_minimum_point_constr}
	A point $x^\star\in{\cal F}\cap {\cal Z} \cap X$ is a local minimum point of Problem (\ref{probconstr}) if, for some $\eps >0$,
	\begin{equation}\label{eq:local_minimum_point_constr}
	\begin{array}{ll}
	f(x^\star)\leq f(x) & \text{ for all } x\in {\cal B}^c(x^\star;\eps)\cap {\cal F}, \\[1.0em]
	f(x^\star)\leq f(x) & \text{ for all } x\in {\cal B}^z(x^\star)\cap {\cal F}.
	\end{array}
	\end{equation}
\end{definition}

Exploiting the necessary optimality conditions introduced in \cite{fasano.2014}, we state the following KKT stationarity definition.

\begin{definition}[KKT Stationary point]\label{def:stationary_point_constr}
	A point $x^* \in {\cal F} \cap {\cal Z} \cap X$ is a KKT stationary point of Problem \eqref{probconstr} 
	if there exists a vector
	$\lambda^\star\in \mathbb{R}^m$ such that, for every $s\in D^c(x^\star)$,
	\begin{eqnarray}
	&&\max \left\{ \xi^\top s\, :\, \xi\in \partial_c f(x^\star) + \sum^m_{i=1}\lambda_i^\star\partial_c g_i(x^\star)\right \}\geq 0,\label{derdir} \\ 
	\bigskip
	&&(\lambda^\star)^Tg(x^\star)=0 \ \text{ and } \ \lambda^\star\ge 0,\label{cmp}
	\end{eqnarray}
	and
	\begin{eqnarray}
	\hspace{-0.30cm}f(x^\star)\le f(x)
	\ \text{ for all }\ x\in {\cal B}^z(x^\star)\cap{\cal F}.\label{eq:local_minimum_point_discr_constr_KKT}\\\nonumber
	\end{eqnarray}
\end{definition}
Similarly to Section \ref{subsec:def_bound_constr}, we can also define weak local minima and weak KKT stationary points by suitably replacing ${\cal B}^z(x^*)$ with ${\cal B}^z(x^*,\tilde D^z)$ in the above two definitions.


\section{An algorithm for bound constrained problems}\label{sec:bound_constrained_case}
In this section, we propose an algorithm for solving the mixed-integer bound constrained problem defined by Problem \eqref{probbox} and we analyze its convergence properties. The optimization over the continuous and discrete variables is performed by means of two local searches based on linesearch algorithms that explore the feasible search directions similarly to the procedures proposed in \cite{fasano.2014,liuzzi.2012,liuzzi.2014,lucidi.2002}. In particular, the \emph{Projected Continuous Search} described in Algorithm~\ref{projContSearch} and the \emph{Discrete Search} described in Algorithm~\ref{DiscrSearch} are the methods adopted to investigate the directions associated with the continuous and discrete variables, respectively. The idea behind the line searches is to return a positive stepsize $\alpha$, namely to update the current iterate, whenever a point providing a sufficient reduction of the objective function is found. In Algorithm~\ref{projContSearch} the sufficient decrease is controlled by the parameter $\alpha$, while in Algorithm~\ref{DiscrSearch} the same role is played by $\xi$. Once such a point is determined, an expansion step is performed in order to explore if the sufficient reduction may be achieved through a larger stepsize.  

\begin{algorithm}[H]
	\caption{Projected Continuous Search ($\tilde\alpha,w,p;\alpha$,$\tilde p$)}\label{projContSearch}
	
	\begin{algorithmic}[1]
		\medskip
		\item[]   {\bf Data.} $\gamma > 0$, $\delta \in (0,1)$.
		
		\medskip
		
		\item[] {\bf Step 0.} Set $\alpha = \tilde\alpha$.
		
		\item[] {\bf Step 1.} {\bf If} $f([w+\alpha p]_{[l,u]})\le f(w)-\gamma\alpha^2$ {\bf
			then} set $\tilde p = p$ go to Step 4.
		
		\item[] {\bf Step 2.} {\bf If} $f([w-\alpha p]_{[l,u]})\le f(w)-\gamma\alpha^2$ {\bf then} set $\tilde p = -p$ and go to Step
		4.
		
		\item[] {\bf Step 3.} Set $\alpha = 0$, {\bf return} $\alpha$ and $\tilde p = p$.
		
		\item[] {\bf Step 4.} Let $\beta={\alpha}/{\delta}$.
		
		\item[] {\bf Step 5.} {\bf If} $f([w+\beta \tilde p]_{[l,u]})>f(w)-\gamma\beta^2$  {\bf return} $\alpha$, $\tilde p$.
		
		\item[] {\bf Step 6.} Set $\alpha=\beta$ and go to
		Step 4.
		
		\par\bigskip\noindent
		
	\end{algorithmic}
\end{algorithm}

\begin{algorithm}[H]
	\caption{Discrete Search ($\tilde\alpha,w,p,\xi;\alpha$)}\label{DiscrSearch}
	
	\begin{algorithmic}[1]
		\medskip
		\item[]   {\bf Data.} $\gamma > 0$.
		
		\medskip
		
		\item[] {\bf Step 0.} Compute the largest $\bar \alpha$  such that $w +\bar \alpha p \in X\cap\mathcal{Z}$.\par
		\item[] \qquad \qquad Set $\alpha =\min\{\bar\alpha,\tilde \alpha \}$.
		\smallskip
		\item[] {\bf Step 1.} {\bf If} $\alpha >0$ and $f(w+\alpha p)\le f(w)-\xi$ {\bf
			then} go to Step 2.\par
		\item[] \qquad \qquad {\bf Else} Set $\alpha = 0$, {\bf return} $\alpha$.
		\smallskip
		\item[] {\bf Step 2.} Set $\beta=\min\{\ \bar \alpha,2{\alpha}\}$.
		\smallskip		
		\item[] {\bf Step 3.} {\bf If} $f(w+\beta p)>f(w)-\xi$,  {\bf return} $\alpha$.
		\smallskip		
		\item[] {\bf Step 4.} Set $\alpha=\beta$ and go to
		Step 2.
		
		\par\bigskip\noindent
		
	\end{algorithmic}
\end{algorithm}

\par\bigskip\bigskip

The algorithm for bound constrained problems proposed in this section performs an alternate minimization along continuous and discrete variables and is hence divided into two phases (see Algorithm~\ref{alg_mixed_3} for a detailed scheme). Starting from a point $x_0 \in X\cap\mathcal{Z}$, in Phase 1 the minimization over the continuous variables is performed by using the Projected Continuous Search (\textbf{Step 7}). If the line search fails, i.e., $\alpha^c_k = 0$, the tentative step for continuous search directions is reduced (\textbf{Step 9}), otherwise the current iterate is updated (\textbf{Step 11}). Then, in Phase 2.A, the directions in the set $D \subset D^z(\tilde x_k)$, where $\tilde x_k$ is the current iterate obtained at the end of Phase 1, are investigated through the Discrete Search (\textbf{Step 16}). If the stepsize returned by the line search performed along a given primitive direction is 0, the corresponding tentative step is halved (\textbf{Step 18}), otherwise the current iterate is updated (\textbf{Step 20}). The directions in $D$ are explored until either a point leading to a sufficient decrease in the objective function is found or $D$ contains no direction to explore. Note that the strategy starts with a subset of $D^z(x_0)$, namely $D$, and  gradually adds directions to it (see Phase 2.B) throughout the iterations. This choice enables the algorithm to reduce the computational cost.
If $\tilde x_k$ obtained at the end of Phase 1 is not modified in Phase 2.A or a direction in $D_k$ along which the Discrete Search does not fail with $\tilde\alpha_{k}^{(d)}=1$ exists, $D_{k+1}$ is set equal to $D_k$ and $D$ is set equal to $D_{k+1}$ (\textbf{Step 32}). Otherwise, $\xi_k$ is reduced (\textbf{Step 24}) and, if all the feasible primitive discrete directions at $\tilde x_k$ have been generated, $D_{k+1}$ and $D$ are not changed compared to the previous iteration (\textbf{Step 26}). Instead, when $D_k \subset D^z(\tilde x_k)$, $D_k$ is enriched with new feasible primitive discrete directions (\textbf{Steps 28--29}) and the initial tentative steps of the new directions are set equal to 1.

We point out that the positive parameter $\xi_k$ plays an important role in the algorithm, since it rules the sufficient decrease of the objective function value within the Discrete Search. The update of the parameter is performed at \textbf{Step 24}. More specifically, we shrink the value of the parameter when both the current iterate is not updated in Phase 2.A and the Discrete Search fails (with  $\tilde\alpha_{k}^{(d)}=1$) along each direction in $D_k$.  
\begin{algorithm}
	\caption{
		DFNDFL
	}\label{alg_mixed_3}
	\begin{algorithmic}[1]
		\par\medskip
		\Algphase{\quad \ \fbox{DATA}}
		\smallskip
		\item \ Let $x_0\in X\cap\mathcal{Z}$, \ $\xi_0>0$, \  $\theta\in (0,1)$;
		\smallskip
		\item \ let $\{s_k\}$ be a sequence such that $s_k\in D^c(x_0)$ and  $\|s_k\|=1$ for all $k$;
		\smallskip
		\item \ let $\tilde\alpha_0^{c} = 1$ be the initial stepsize along $s_k$;
		\smallskip
		\item \ let $D=D_0\subset D^z(x_0)$ be a set of initial feasible primitive discrete directions at point $x_0$;
		\smallskip
		\item \ let $\tilde\alpha_0^{(d)} = 1$ be the initial stepsizes along  $d \in D$.
		
		\par\medskip
		
		\item {\bf For} $k=0,1,\dots$
		\par\medskip
		
		\Algphase{\qquad\quad\fbox{PHASE 1 - Explore continuous variables }}
		\smallskip
		\item \qquad Compute $\alpha_k^{c}$ and $\tilde s_k$ by the {\em Projected Continuous Search}$(\tilde\alpha_k^{c},x_k,s_k;\alpha_k^{c}, \tilde s_k)$.
		\smallskip 
		\item \qquad {\bf If} $(\alpha_k^{c} = 0)$ {\bf then}		
		\item \qquad\qquad $\tilde \alpha_{k+1}^{c}=\theta \tilde \alpha_k^{c}$ and $\tilde x_k= x_k$,
		\item \qquad {\bf else} 
		\item \qquad\qquad $\tilde\alpha_{k+1}^{c}= \alpha_k^{c}$
		and $\tilde x_k=[x_k+\alpha_k^{c} \tilde s_k]_{[l,u]}$. 
		\smallskip
		\item \qquad {\bf End If} 		
		\medskip
		\Algphase{\qquad\quad\fbox{PHASE 2.A - Explore discrete variables}}
		\smallskip
		\item \qquad Set $y^+=\tilde x_k$. 
		
		\par\medskip
		\item \qquad {\bf While} $D\neq \emptyset$ and $y^+=\tilde x_k$ {\bf do}
		\medskip
		
		
		\item \qquad\quad Choose $d \in D$, set $D=D\setminus\{d\}$ and $y=y^+$. 
		
		\smallskip
		
		\item \qquad\quad Compute $\alpha$ by the {\em Discrete Search}$(\tilde\alpha_k^{(d)},y,d,\xi_k;\alpha)$.
		
		\medskip
		%
		
		\item\label{alg1:line7} \qquad\quad {\bf If} $\alpha = 0$ {\bf then} 
		\item \qquad\qquad\quad Set $y^+=y$ \ \ and\ \   $\tilde\alpha_{k+1}^{(d)} =\max\{1,\lfloor\tilde\alpha_k^{(d)}/2\rfloor\}$,
		\item \qquad\quad {\bf else}
		\item \qquad\qquad\quad Set $y^+=y+\alpha d$ 
		\ \ and $\tilde\alpha_{k+1}^{(d)} = \alpha$.
		\smallskip
		
		\par
		
		\item \qquad\quad {\bf End If}
		
		\par\medskip
		
		\item \qquad {\bf End While}
		\medskip
		
		\Algphase{\qquad\quad\fbox{PHASE 2.B - Update the set of discrete search directions}}
		\smallskip
		\item \qquad {\bf If} $y^+ = \tilde x_k$  and {\em  Discrete Search} fails with $\tilde\alpha_{k}^{(d)} = 1$ for all $d\in D_{k}$ \textbf{then}
		
		\item \qquad\quad  Set $\xi_{k+1}=\theta \xi_k$. 
		
		\item \qquad\quad {\bf If} $ D_{k}\supseteq D^z(\tilde x_k)$ {\bf then}
		
		\item \qquad\qquad\quad      Set $D_{k+1} = D_k$ and  $D=D_{k+1}$,

		\item \qquad\quad {\bf else}
		
		\item \qquad\qquad\quad Generate
		$D_{k+1}$  such that $D_{k+1} \subseteq D^z(\tilde x_k)$ and $ D_{k+1} \supset D_k$, set
		$D=D_{k+1}$.
		
		\item \qquad\qquad\quad Set $\tilde\alpha_{k+1}^{(d)} = 1$  for all $d\in D_{k+1}\setminus D_k$.		
		
		\item \qquad\quad {\bf End If}
		
		\item \qquad {\bf else}
		\item \qquad\quad  Set $D_{k+1} = D_k$ and $D=D_{k+1}$.
		
		\item \qquad {\bf End If}
		\medskip
		\Algphase{\qquad\quad\fbox{PHASE 3 - Update iterates}}
		\smallskip
		\item \qquad Find $x_{k+1}\in X \cap {\cal Z}$ such that $f(x_{k+1})\le f(y^+)$.
		\medskip
		\item {\bf End For}
		
		\par\bigskip\noindent
		
	\end{algorithmic}
\end{algorithm}


The following propositions guarantee that the algorithm is well-defined.
\begin{proposition}\label{prop:infinite_cycle_alg_1}
	The {\em Projected Continuous Search} cannot infinitely cycle between Step 4 and Step 6.
\end{proposition}
\proof \quad We assume by contradiction that in the Projected Continuous Search an infinite monotonically increasing sequence of positive numbers $\{\beta_j\}$ exists such that $\beta_j \to \infty$ for $j \to \infty$ and
\[
f([w+\beta_j p]_{[l,u]})\leq f(w)-\gamma\beta^2_j.
\]
Since by the instructions of the procedures we have that $[w+\beta_j p]_{[l,u]} \in X \cap {\cal Z}$, the previous relation is in contrast with the compactness of $X$, by definition of compact set, and with the continuity of function $f$. These arguments conclude the proof. 
$\hfill\Box$ \\
\par
\begin{proposition}\label{prop:infinite_cycle_alg_2}
	The {\em Discrete Search} cannot infinitely cycle between Step 2 and Step 4.
\end{proposition}
\proof \quad We assume by contradiction that in the Discrete Search an infinite monotonically increasing sequence of positive numbers $\{\beta_j\}$ exists such that $\beta_j \to \infty$ for $j \to \infty$ and
\[
f(w+\beta_j p)\leq f(w)-\xi.
\]
Since by the instructions of the procedures we have that $w+\beta_j p \in X \cap {\cal Z}$, the previous relation is in contrast with the compactness of $X$, by definition of compact set. This argument concludes the proof. 
$\hfill\Box$ \\
\par
In the following proposition, we prove that the Projected Continuous Search returns stepsizes that eventually go to zero.

\begin{proposition}\label{prop:maxtozerosi}
	Let $\{\alpha_k^{c}\}$ and $\{\tilde\alpha_k^{c}\}$ be the sequences yielded by Algorithm DFNDFL for each $s_k\in D^c(x_k)$. Then
	\begin{equation}\label{maxtozerosi}
		\lim_{k\to\infty}\max \{\alpha_k^{c}, \tilde \alpha_k^{c}\} = 0.
	\end{equation}
\end{proposition}
\proof  \quad The proof follows with minor modifications from the proof of Proposition 2.5 in \cite{fasano.2014}.

$\hfill\Box$ 

\par\medskip
\begin{proposition}\label{tech:lemma1}
	Let $\{\xi_k\}$ be the sequence produced by
	Algorithm DFNDFL. Then	
	\begin{equation*}
		\lim_{k\to\infty} \xi_k = 0.
	\end{equation*}
\end{proposition}
\proof \quad 
By the instruction of Algorithm DFNDFL, it follows that
$0< \xi_{k+1} \leq \xi_k$ for all $k$, meaning that the sequence $\{\xi_k\}$ is monotonically nonincreasing. Hence, $\{\xi_k\}$ converges to a limit $M \geq 0$. Suppose, by contradiction, that
$M > 0$. This implies that an index $\bar k > 0$ exists such that $\xi_{k+1} = \xi_k = M$ for all $k\geq\bar k$.
Moreover, for every index $k\geq\bar k$, a direction $d \in D^z(\tilde x_k)$ exists such that
\begin{equation}\label{DFL_eqassurdo}
	f(x_{k+1}) \leq f(\tilde x_k + \alpha_k^{(d)}d) \leq f(\tilde x_k) - M \leq f(x_k) - M,
\end{equation}
In fact, if such an index did not exist, the algorithm would set $\xi_{k+1} = \theta \xi_k$. Relation (\ref{DFL_eqassurdo}) implies $f(x_k)\to -\infty$, which is in contrast with the assumption that $f$ is continuous on the compact set $X$. This concludes the proof. $\hfill\Box$
\par\medskip
\begin{remark}\label{remark2} By the preceding proposition and the updating rule of the parameter $\xi_k$ used in Algorithm DFNDFL, it follows that the set $$H = \{k:\xi_{k+1} < \xi_k\}$$ is infinite.  
\end{remark} 
\par
The previous result is used to prove the next lemma, which in turn is essential to prove the global convergence result related to the continuous variables. This lemma states that the asymptotic convergence properties of the sequence $\{s_k\}$ still hold when the projection operator is adopted. Its proof closely resembles the proof in \cite[Lemma 2.6]{fasano.2014}. 
\begin{lemma}\label{lemmaseq}
	Let $\{x_k\}$ and $\{s_k\}$ be the sequence of points and the sequence of continuous search directions produced by Algorithm DFNDFL, respectively, and $\{\eta_k\}$ be a sequence such that $\eta_k>0$, for all $k$. Further, let $K$ be a subset of indices such that
	\begin{eqnarray}
		\lim_{k \to \infty, k \in K} x_k = \bar x, && \label{limxk}\\
		\lim_{k \to \infty, k \in K} s_k = \bar s, && \label{limdk2}\\
		\lim_{k \to \infty, k \in K} \eta_k = 0. && \label{limalpha1}
	\end{eqnarray}
	with $\bar x\in X\cap\mathcal{Z}$ and $\bar s \in D^c(\bar x)$, $\bar s\neq 0$. Then,
	\begin{itemize}
		\item[(i)] for all $k \in K$ sufficiently large,
		
		$$[x_k+\eta_k s_k]_{[l,u]}\neq x_k,$$
		
		\item[(ii)] the following limit holds
		$$\displaystyle\lim_{k\to\infty, k\in K}v_k = \bar s,$$
		where
		\begin{equation}\label{vk}
			v_k=\frac{[x_k+\eta_k s_k]_{[l,u]}-x_k}{\eta_k}.
		\end{equation}
		
	\end{itemize}

\end{lemma}
\proof \ The proof can be obtained by suitably adapting the proof of Lemma 2.6 in \cite{fasano.2014}. $\hfill\Box$ 

\par\medskip

The main convergence result related to the continuous variables is proved in the next proposition. It basically states that every limit point of the subsequence of iterates defined by the set $H$ (see Remark~\ref{remark2}) is a stationary point with respect to the continuous variables.
\begin{proposition}\label{tech:lemma3-3}
	Let $\{x_k\}$ be the sequence of points produced by Algorithm DFNDFL. Let $H\subseteq\{1,2,\dots\}$ be defined as in Remark
	\ref{remark2} and let $\bar x \in X \cap {\cal Z}$ be any accumulation point of $\{x_k\}_H$. If the subsequence $\{s_k\}_H$, with $(s_k)_i=0$ for $i \in I^z$, is dense in the unit sphere (see Definition \ref{def:dense_sequence}), then $\bar x$ satisfies
	\begin{align}
		&f^{\circ}_c(\bar x; s)\geq 0\quad\quad \mbox{for all}\ s \in D^c(\bar x).\label{nulgrad_unc2-3}
	\end{align}
\end{proposition}
\proof \
For any accumulation point $\bar x$ of $\{x_k\}_H$, let $K\subseteq H$ be an index set such that
\begin{equation}\label{limK2-3}
	\lim_{k\to\infty,k\in K}x_k = \bar x.
\end{equation}
Notice that, for all $k\in K$, $(\tilde x_k)_z = (x_k)_z$ and $\tilde\alpha_k^{(d)}=1$, $d \in D_k$, by the instructions of Algorithm DFNDFL. Hence, for all $k\in K$, by recalling (\ref{limK2-3}), the discrete variables are no longer updated. 

Now, reasoning as in the proof of \cite[Proposition 2.7]{fasano.2014} one can show that no direction $\bar s\in D^c(\bar x)\cap S(0,1)$ can exist such that
\begin{equation}\label{contradict1_2}
	f^\circ_c(\bar x;\bar s) < 0.
\end{equation}		
$\hfill\Box$

\par\medskip

The next proposition states that every limit point of the subsequence of iterates defined by the set $H$ (see Remark~\ref{remark2}) is a local minimum with respect to the discrete variables.
\begin{proposition}\label{tech:lemma2-3}
	Let $\{x_k\}$, $\{\tilde x_k\}$, and $\{\xi_k\}$ be the sequences produced by Algorithm DFNDFL. Let
	$H\subseteq\{1,2,\dots\}$ be defined as in Remark \ref{remark2} and $x^* \in X \cap {\cal Z}$ be any accumulation point of $\{x_k\}_H$, then
	\[
	f(x^*)\leq f(\bar x),\qquad \mbox{for all}\ \bar x \in {\cal B}^z(x^*).
	\]
\end{proposition}
\proof \quad Let $K \subseteq H$ be an index set such that
\[
\lim_{k\to\infty,k\in K}x_k = x^*.
\]
For every $k\in K\subseteq H$, we have
\begin{eqnarray*}
	(\tilde x_k)_z  & = & (x_k)_z,\\
	\nonumber\tilde\alpha_k^{(d)} & = & 1,\quad d \in D_k,
\end{eqnarray*}
meaning that the discrete variables are no longer updated by the Discrete Search.

Let us consider any point $\bar x\in {\cal B}^z(x^*)$. By the definition of discrete neighborhood ${\cal B}^z(x^*)$, a direction $\bar d \in D^z(x^*)$ exists such that
\begin{equation}\label{DFL_discr1-3}
\bar x = x^* + \bar d.
\end{equation}
Recalling the steps in Algorithm DFNDFL, we have,
for all $k\in H$ and sufficiently large, that
$$
(x^*)_z = (x_k)_z = (\tilde x_k)_z.
$$
Further, by Proposition \ref{prop:maxtozerosi}, we have
$$
\lim_{k\to\infty,k\in K} \tilde x_k = x^*.
$$
Then, for all $k\in K$ and sufficiently large, \eqref{DFL_discr1-3} implies
\[
(x_k + \bar d)_z = (\tilde x_k + \bar d)_z = (x^*  + \bar d)_z = (\bar x)_z.
\]
Hence, for all $k\in K$ and sufficiently large, by the definition of discrete neighborhood we have $\bar d \in D^z(\tilde x_k)$ and
$$
\tilde x_k + \bar d \in X\cap  \mathcal{Z}.
$$
Then, since $k \in K\subseteq H$, by the definition of $H$ we have
\begin{equation}\label{DFL_fallimento_disc-3}
f(\tilde x_k + \bar d) > f(\tilde x_k) -\xi_k.
\end{equation}
Now, by Proposition \ref{tech:lemma1}, and taking the limit for $k\to\infty$, with $k\in K$, in \eqref{DFL_fallimento_disc-3}, the
result follows. $\hfill\Box$ \\
\par
Now we can prove the main convergence result of the algorithm.
\begin{theorem}\label{theo4}
	Let $\{x_k\}$ be the sequence of points generated by Algorithm DFNDFL. Let $H\subseteq\{1,2,\dots\}$ be defined as in Remark
	\ref{remark2} and let $\{s_k\}_H$, with $(s_k)_i=0$ for $i \in I^z$, be a dense subsequence in the unit sphere. Then,
	\begin{itemize}
		\item[(i)] a limit point of $\{x_k\}_H$ exists;
		\item[(ii)] every limit point $x^*$ of $\{x_k\}_H$ is stationary for Problem \eqref{probbox}.
	\end{itemize}
\end{theorem}
\proof \quad As regards point (i), since $\{x_k\}_H$ belongs to the compact set $X \cap {\cal Z}$, it admits limit points. The prove of point (ii)
follows by considering Propositions \ref{tech:lemma3-3} and \ref{tech:lemma2-3}. $\hfill\Box$

\section{An algorithm for nonsmooth nonlinearly constrained problems}\label{sec:nonlinearly_constrained_case}
In this section, we consider the nonsmooth nonlinearly constrained problem defined in Problem~\eqref{probconstr}. The nonlinear constraints are handled through a simple penalty approach (see, e.g., \cite{fasano.2014}). In particular, given a positive parameter $\varepsilon > 0$, we introduce the following penalty function
\[
P(x;\varepsilon) = f(x) + \frac{1}{\varepsilon}\sum_{i=1}^m\max\left\{0,g_i(x)\right\},
\]
which allows us to define the following bound constrained problem
\begin{equation}\label{prob_Z}
\begin{array}{ll}
\min & P(x;\varepsilon)\\
s.t. & x\in X \cap \mathcal{Z}.
\end{array}
\end{equation}
Hence, only the nonlinear constraints are penalized and the minimization is performed over the set $X \cap \mathcal{Z}$. The algorithm described in Section~\ref{sec:bound_constrained_case} is thus suited for solving this problem, as highlighted in the following remark.
\begin{remark}
	Observe that, for any $\varepsilon>0$, the structure and properties of Problem (\ref{prob_Z}) are the same as Problem (\ref{probbox}).	
	The Lipschiz continuity (with respect to the continuous variables) of the penalty function $P(x;\varepsilon)$ follows by the Lipschitz continuity of $f$ and $g_i$, with $i \in \{1,\dots,m\}$. In particular, called $L_f$ and $L_{g_i}$ the Lipschitz constants of $f$ and $g_i$, respectively, we have that the Lipschitz constant of the penalty function $P(x;\varepsilon)$ is
	\[
	L \leq L_f + \frac{1}{\varepsilon}\sum_{i=1}^m L_{g_i}.
	\]
\end{remark}
To prove the equivalence between Problem \eqref{probconstr} and Problem \eqref{prob_Z}, we report an extended version of the  Mangasarian-Fromowitz Constraint Qualification (EMFCQ) condition for Problem \eqref{probconstr}, which takes into account its mixed-integer structure. This condition states that at a point that is infeasible for Problem~\eqref{probconstr}, a direction feasible with respect to $X \cap \mathcal{Z}$ (according to Definitions \ref{def:set_feasible_primitive_discrete_directions} and \ref{def:cone_feasible_continuous_directions}) that guarantees a reduction in the constraint violation exists.
%

\begin{assumption}[EMFCQ for mixed-integer problems]\label{assmfcq}
	Given Problem \eqref{probconstr}, for any \\ $x\in (X\cap\mathcal{Z})\setminus\stackrel{\circ}{\cal F}$, one of the following conditions holds:
	\begin{itemize}
		\item[(i)] a direction $s\in D^c(x)$ exists such that  \[
		(\xi^{g_i})^\top s < 0,
		\]
		for all $\xi^{g_i}\in\partial_c g_i(x)$ with $i\in\{h \in \{1,\dots,m\}: \ g_h(x)\geq 0\}$; 
		\item[(ii)] a direction $\bar d \in D^z(x)$ exists such that
		$$\sum_{i=1}^m \max\{0, g_i(x+\bar d)\} < \sum_{i=1}^m \max\{0, g_i(x)\}.$$ 
	\end{itemize}
\end{assumption}


In order to prove the main convergence properties of the algorithm in this case, we first need to establish the equivalence between the original constrained Problem \eqref{probconstr} and the penalized Problem \eqref{prob_Z}. The proof of this result is very technical and quite similar to analogous results from \cite{fasano.2014}. We report it in the Appendix for the sake of major clarity.

\par

Exploiting this technical result, we can apply the algorithm proposed in Section \ref{sec:bound_constrained_case} to solve Problem~\eqref{prob_Z}, provided that the penalty parameter is sufficiently small, as stated in the next proposition. The algorithmic scheme designed for solving Problem \eqref{probconstr} is obtained from Algorithm~DFNDFL by replacing $f(x)$ with $P(x;\varepsilon)$, where $\varepsilon>0$ is a sufficiently small value. We point out that in this new scheme  both the linesearch procedures are  performed by replacing $f(x)$ with $P(x;\varepsilon)$ as well. We refer to this new scheme as DFNDFL--CON. 
\begin{proposition}\label{mainconv_constr} Let Assumption \ref{assmfcq} hold and let $\{x_k\}$ be the sequence produced by Algorithm~DFNDFL--CON.
Let $H\subseteq\{1,2,\dots\}$ be defined as in Remark
	\ref{remark2} and let $\{s_k\}_H$, with $(s_k)_i=0$ for $i \in I^z$, be a dense subsequence in the unit sphere. Then,
$\{x_k\}_H$ admits limit points. Furthermore, a threshold value $\varepsilon^*$ exists such that
for all $\varepsilon \in (0, \varepsilon^*]$ every limit point $x^*$ of $\{x_k\}_H$ is stationary for Problem \eqref{probconstr}.
\end{proposition}
\proof \quad The proof follows from Proposition \ref{equivalenza} and Theorem \ref{theo4}. $\hfill\Box$ \\
\par

\section{Numerical experiments}\label{sec:numerical_exp}
In this section, we report the results of the numerical experiments performed on a set of test problems selected from the literature. 
In particular, state-of-the-art solvers are used as benchmarks to test the efficiency and reliability of the proposed algorithm. First we consider the bound constrained case, then we move on to nonlinearly constrained problems. In both cases, to improve the performance of DFNDFL, a modification to Phase 1 is introduced by drawing inspiration from the algorithm CS-DFN proposed in~\cite{fasano.2014} for continuous nonsmooth problems. In particular, recalling that $I^c \cup I^z = \{1,\ldots,n\}$, the change consists in investigating the set of coordinate directions $\{\pm e^1, \ldots, \pm e^{|I^c|}\}$ before exploring a direction from the sequence $\{s_k\}$. Since this set is constant over the iterations, the actual and tentative stepsizes $\alpha_k^{(i)}$ and $\tilde\alpha_k^{(i)}$ can be stored for each coordinate direction $i$, with $i \in \{1,\ldots,n\}$. These stepsizes are reduced whenever the continuous line search (i.e., Algorithm \ref{projContSearch} without the projection operator) does not determine any point that satisfies the sufficient decrease condition. When their values become sufficiently small, a direction from the dense sequence $s_k$ is explored. This improvement allows the algorithm to benefit from the presence of the stepsizes $\alpha_k^{(i)}$ and $\tilde\alpha_k^{(i)}$, whose values depend on the knowledge across the iterations of the sensitivity of the objective function over the coordinate directions. Therefore, the efficiency of the modified DFNDFL is expected to be higher. 
\par
The codes related to the DFNDFL and DFNDFL--CON algorithms, together with the test problems used in the experiments are freely available for download at the \emph{Derivative-Free Library} web page
\begin{center}
	\url{http://www.iasi.cnr.it/~liuzzi/DFL/}
\end{center}

\subsection{Algorithms for benchmarking}
The algorithms selected as benchmarks  are listed below:
\begin{itemize}
	\item DFL (see \cite{liuzzi.2012}), a derivative-free linesearch algorithm for bound constrained (in this case the algorithm is referred to as DFL {\em box}) and nonlinearly constrained problems (referred to as DFL {\em gen}).
	\item RBFOpt (see \cite{costa.2018}), an open-source library RBFOpt for solving black-box bound constrained optimization problems with expensive function evaluations.
	\item NOMAD v.3.9.1 (see \cite{Nomad}), a software package which implements the mesh adaptive direct search algorithm. 
	\item MISO (Mixed-Integer Surrogate Optimization) framework, a model-based approach using surrogates \cite{muller.2015}.
	
\end{itemize} 
All the algorithms reported above support mixed-integer problems, thus being suited for the comparison with the algorithm proposed in this work. The maximum allowable number of function evaluations in each experiment is 5000. All the codes have been run on an Intel Core i7 10th generation CPU PC running Windows 10 with 16GB of memory. More precisely, all the test problems, DFNDFL, DFL {\em box}, DFL {\em gen} and RBFOpt are coded in python and have been run using python 3.8. NOMAD, on the other hand, is delivered as a collection of C++ codes and has been run using the provided PyNomad python interface. As for MISO, it is coded in matlab and has been run using Matlab R2020b but using the python coded problems through the matlab python engine.
\par
As regards the parameters used in both DFNDFL and DFL, the values used in the experiments are $\gamma = 10^{-6}, \, \delta = 0.5, \, \xi_0 = 1$, and $ \theta = 0.5$. Moreover, the initial tentative steps along the coordinate directions $\pm e^i$ and $s_k$ of the modified DFNDFL are 
\begin{eqnarray*}
	\tilde\alpha_0^i & = & (u^i - \ell^i)/2 \qquad \text{ for all } i \in I^c,\\
	\tilde\alpha_0   & = & \frac{1}{n}\sum_{i=1}^n\tilde\alpha_0^i, 
\end{eqnarray*}
while for the discrete directions $d$ the initial tentative stepsize $\tilde\alpha_0^{(d)}$ is fixed to 1. 

Another computational aspect that needs to be further discussed is the generation of the continuous and discrete directions. Indeed, in Phases 1 and 2 of DFNDFL, new search directions might be generated to thoroughly explore neighborhoods of the current iterate. To this end, a dense sequence of directions $\{s_k\}$ is required in Phase 1 to explore the continuous variables. Similarly, in Phase 2, new primitive discrete directions must be generated when some suitable conditions hold. In both the situations, we used Sobol sequences (\cite{bratley.1988,sobol.1976}).

\par
As concerns the parameters used for running RBFOpt and NOMAD, while the former is executed by using the default values, for the latter two different algorithms are considered. The first one is based on the default settings, while the second one results from disabling the usage of models in the search phase, which precisely is obtained by setting {\tt DISABLE MODELS}. This second version is denoted in the remainder of this document as NOMAD (w/o mod).

\subsection{Test problems}
The comparison between Algorithm~DFNDFL and some state-of-the-art solvers is reported for 24 bound constrained problems. The first 16 problems, which are related to minimax and nonsmooth unconstrained optimization problems, have been selected from \cite[Sections 2 and 3]{luksan.2000}, while the remaining 8 problems have been chosen from \cite{muller.2015,muller.2013b}. The problems are listed in Table~\ref{tab:boxprob} along with the respective number of continuous ($n_c$) and discrete ($n_z$) variables.

\begin{table}[th]
	\begin{center}
		\begin{tabular}{lcrrlcrr}\hline
			Problem name & Source & $n_c$ & $n_z$ & Problem name & Source & $n_c$ & $n_z$ \\\hline
			gill              & \cite{luksan.2000} & 5 & 5 & 
			goffin            & \cite{luksan.2000} &25 &25 \\
			l1hilb(20)        & \cite{luksan.2000} &10 &10 & 
			l1hilb(30)        & \cite{luksan.2000} &15 &15 \\
			l1hilb(40)        & \cite{luksan.2000} &20 &20 & 
			l1hilb(50)        & \cite{luksan.2000} &25 &25 \\
			maxl              & \cite{luksan.2000} &10 &10 & 
			maxq(20)          & \cite{luksan.2000} &10 &10 \\
			maxq(30)          & \cite{luksan.2000} &15 &15 & 
			maxq(40)          & \cite{luksan.2000} &20 &20 \\
			maxq(50)          & \cite{luksan.2000} &25 &25 & 
			maxquad           & \cite{luksan.2000} & 5 & 5 \\
			mxhilb            & \cite{luksan.2000} &25 &25 & 
			osborne 2         & \cite{luksan.2000} & 6 & 5 \\
			polak 2           & \cite{luksan.2000} & 5 & 5 & 
			polak 3           & \cite{luksan.2000} & 6 & 5 \\
			shell dual        & \cite{luksan.2000} & 8 & 7 & 
			steiner 2         & \cite{luksan.2000} & 6 & 6 \\
			tr48              & \cite{luksan.2000} &24 &24 & 
			watson            & \cite{luksan.2000} &10 &10 \\
			wong2             & \cite{luksan.2000} & 5 & 5 & 
			wong3             & \cite{luksan.2000} &10 &10 \\
			SO-I prob. 7      & \cite{muller.2013b}& 5 & 5 & 
			SO-I prob. 9      & \cite{muller.2013b}& 6 & 6 \\
			SO-I prob.10      & \cite{muller.2013b}&15 &15 & 
			SO-I prob.13      & \cite{muller.2013b}& 5 & 5 \\
			SO-I prob.15      & \cite{muller.2013b}& 6 & 6 & 
			MISO prob. 6      & \cite{muller.2015} & 8 & 7 \\
			MISO prob. 8      & \cite{muller.2015} & 5 &10 & 
			MISO prob.10      & \cite{muller.2015} &30 &30 \\
			\cline{1-8}
		\end{tabular}
	\end{center}
	\caption{Bound constrained test problems collection}\label{tab:boxprob}
\end{table}
Since the problems from \cite{luksan.2000} are unconstrained, in order to suit the class of problems addressed in this work, the following bound constraints are considered for each variable
\[
\ell^i = (\tilde x_0)^i -10 \leq \tilde x^i \leq (\tilde x_0)^i + 10 = u^i \quad \text{ for all } i \in \{1,\dots,n\},
\]
where $\tilde x_0$ is the starting point. Furthermore, since the problems in \cite{luksan.2000} have only continuous variables, the rule applied to obtain mixed-integer problems is to consider a number of integer variables equal to $n_z = \left\lfloor{n/2}\right\rfloor$ and a number of continuous variables equal to $n_c = \left\lceil{n/2}\right\rceil$, where $n$ denotes the dimension of each original problem and $\left\lfloor{\cdot}\right\rfloor$ and $\left\lceil{\cdot}\right\rceil$ are the floor and ceil operators, respectively. 

More specifically, let us consider both the continuous bound constrained optimization problems from~\cite{luksan.2000}, whose formulation is 
\begin{equation}\label{prob:original_cont_bound_constr_probl}
\begin{array}{l}
\min\ \tilde f(\tilde x)\\
s.t.\ \ell^i\leq \tilde x^i\leq u^i \qquad \ \text{ for all }\ i\in \{1,\ldots,n\}, \\
\phantom{s.t.}\ \tilde x_i \in \mathbb{R} \qquad\qquad \ \ \text{ for all }\ i\in \{1,\ldots,n\}, \\
\end{array}
\end{equation}
and the original mixed-integer bound constrained problems from \cite{muller.2013b,muller.2015}, which can be stated as
\begin{equation}\label{prob:original_MINLP_bound_constr_probl}
	\begin{array}{l}
		\min\ \tilde f(\tilde x)\\
		s.t.\ \ell^i\leq \tilde x^i\leq u^i \qquad \ \text{ for all }\ i\in \{1,\ldots,n\}, \\
		\phantom{s.t.}\ x_i \in \mathbb{R} \qquad\qquad\ \ \text{ for all }\ i \in I^c, \\
		\phantom{s.t.}\ x_i \in \mathbb{Z} \qquad\qquad\ \ \text{ for all }\ i \in I^z. \\
	\end{array}
\end{equation}
The resulting mixed-integer problem we deal with in here can be formulated as follows
\begin{equation}\label{prob:MINLP_experimental_results}
\begin{array}{l}
	\min\ f(x)\\
	s.t.\ \ell^i\leq x^i\leq u^i \qquad \ \text{ for all }\ i \in I^c, \\
	\phantom{s.t.}\ 0 \leq x^i \leq 100 \qquad \text{ for all }\ i\in I^z,\\
	\phantom{s.t.}\ x_i \in \mathbb{R} \qquad\qquad\ \ \text{ for all }\ i \in I^c, \\
	\phantom{s.t.}\ x_i \in \mathbb{Z} \qquad\qquad\ \ \text{ for all }\ i \in I^z,
\end{array}
\end{equation}
where $f(x) = \tilde f(\tilde x)$ with 
	\[
	\tilde x^i  = \left\{\begin{array}{ll}
	x^i, &\quad \mbox{for all}\ i \in I^c, \\
	\ell^i + x^i(u^i-\ell^i)/100,                          &\quad \mbox{for all}\ i    \in I^z.\\
	\end{array}\right.
	\]
Moreover, the starting point $x_0$ adopted for Problem (\ref{prob:MINLP_experimental_results}) is
	\[
	(x_0)^i = \left\{\begin{array}{ll}
	(u^i - \ell^i)/2& \qquad \text{ for all } i\in I^c,\\
	50& \qquad \text{ for all } i \in I^z.
	\end{array}\right.
	\]

As concerns constrained mixed-integer optimization problems, 
the performances of Algorithm~DFNDFL--CON are assessed on 204 problems with general nonlinear constraints. Such problems are obtained by adding to the 34 bound-constrained problems defined by Problem (\ref{prob:MINLP_experimental_results}) the 6 classes of constraints reported below   
\begin{center}\small
	\begin{tabular}{ll}
		$g_j(x) = (3-2x_{j+1})x_{j+1} - x_j -2x_{j+2} + 1 \leq 0,$ & for all \ $j \in \{1,\dots,n-2\}$,\\ \smallskip
		$g_j(x) = (3-2x_{j+1})x_{j+1} - x_j -2x_{j+2} + 2.5 \leq 0$, & for all \ $j \in \{1,\dots,n-2\}$,\\ \smallskip
		$ g_j( x) = x_j^2 + x_{j+1}^2 + x_jx_{j+1} -2x_j -2x_{j+1} +1 \leq 0$, & for all \ $j \in \{1,\dots,n-1\}$,\\ \smallskip
		$g_j(x) = x_j^2 + x_{j+1}^2 + x_jx_{j+1} -1 \leq 0$, & for all \ $j \in \{1,\dots,n-1\}$,\\ \smallskip
		$g_j( x) = (3-0.5x_{j+1})x_{j+1} - x_j -2x_{j+2} + 1 \leq 0$, & for all \ $j \in \{1,\dots,n-2\}$,\\ \smallskip
		$g_1( x) = \sum_{j=1}^{n-2}((3-0.5x_{j+1})x_{j+1} - x_j -2x_{j+2} + 1) \leq 0$. & $\phantom{j=\{1,\dots,n-2\}}$.
	\end{tabular}
\end{center}
Thus, the number of general nonlinear constraints ranges from 1 to 59.

\subsection{Data and performance profiles}
The comparison among the algorithms is carried out by using data and performance profiles, which are benchmarking tools widely used in derivative-free optimization (see \cite{more.2009}). In particular, given a set $S$ of algorithms, a set $P$ of problems, and a convergence test, data and performance profiles provide complementary information to assess the relative performance among the different algorithms in $S$ when applied to solve problems in $P$. Specifically, data profiles allow gaining insight on the percentage of problems that are solved (according to the convergence test) by each algorithm within a given budget of function evaluations, while the performance profiles allow assessing how well an algorithm performs with respect to the other. For each $s\in S$ and $p \in P$, the number of function evaluations required by algorithm $s$ to satisfy the convergence condition on problem $p$ is denoted as $t_{p,s}$. Given a tolerance $0< \tau < 1$ and denoted as $f_L$ the smallest objective function value computed by any algorithm on problem $p$ within a given number of function evaluations, the convergence test is 
\[
f(x_k) \leq f_L + \tau(\hat f(x_0) - f_L),
\]
where $\hat f(x_0)$ is the objective function value of the worst feasible point determined by all
the solvers (note that in the bound-constrained case, $\hat f(x_0) = f (x_0)$). The above convergence test 
requires the best point to achieve a sufficient reduction from the value $\hat f(x_0)$ of the objective function at the starting point. When dealing with problems with general constraints, we set to $+\infty$ the value of the objective function at infeasible points. Note that the smaller the value of the tolerance $\tau$ is, the higher accuracy is required at the best point. In particular, three levels of accuracy are considered in this paper for the parameter $\tau$, namely, $\tau \in \{10^{-1}, 10^{-3},10^{-5}\}$.

Performance and data profiles of solver $s$ can be formally defined as follows
\begin{eqnarray*}
	\rho_s(\alpha) & = & \frac{1}{|P|}\left|\left\{p\in P: \frac{t_{p,s}}{\min\{t_{p,s'}:s'\in S\}}\leq\alpha\right\}\right|,\\
	d_s(\kappa) & = & \frac{1}{|P|}\left|\left\{p\in P: t_{p,s}\leq\kappa(n_p+1)\right\}\right|,
\end{eqnarray*}
where $n_p$ is the dimension of problem $p$. While $\alpha$ indicates that the number of function evaluations required by algorithm $s$ to achieve the best solution is $\alpha$--times the number of function evaluations needed by the best algorithm, 
$k$ denotes the number of simplex gradient estimates, with $n_p + 1$ being the number associated with one simplex gradient. Important features for the comparison are $\rho_s(1)$, which is a measure of the efficiency of the algorithm, since it is the percentage of problems for which the algorithm $s$ performs the best, and the height reached by each profile as the value of $\alpha$ or $k$ increases, which measures the reliability of the solver.

\subsubsection*{The bound constrained case}

\begin{figure}[htbp]
	\centering
	\includegraphics[width=0.80\textwidth]{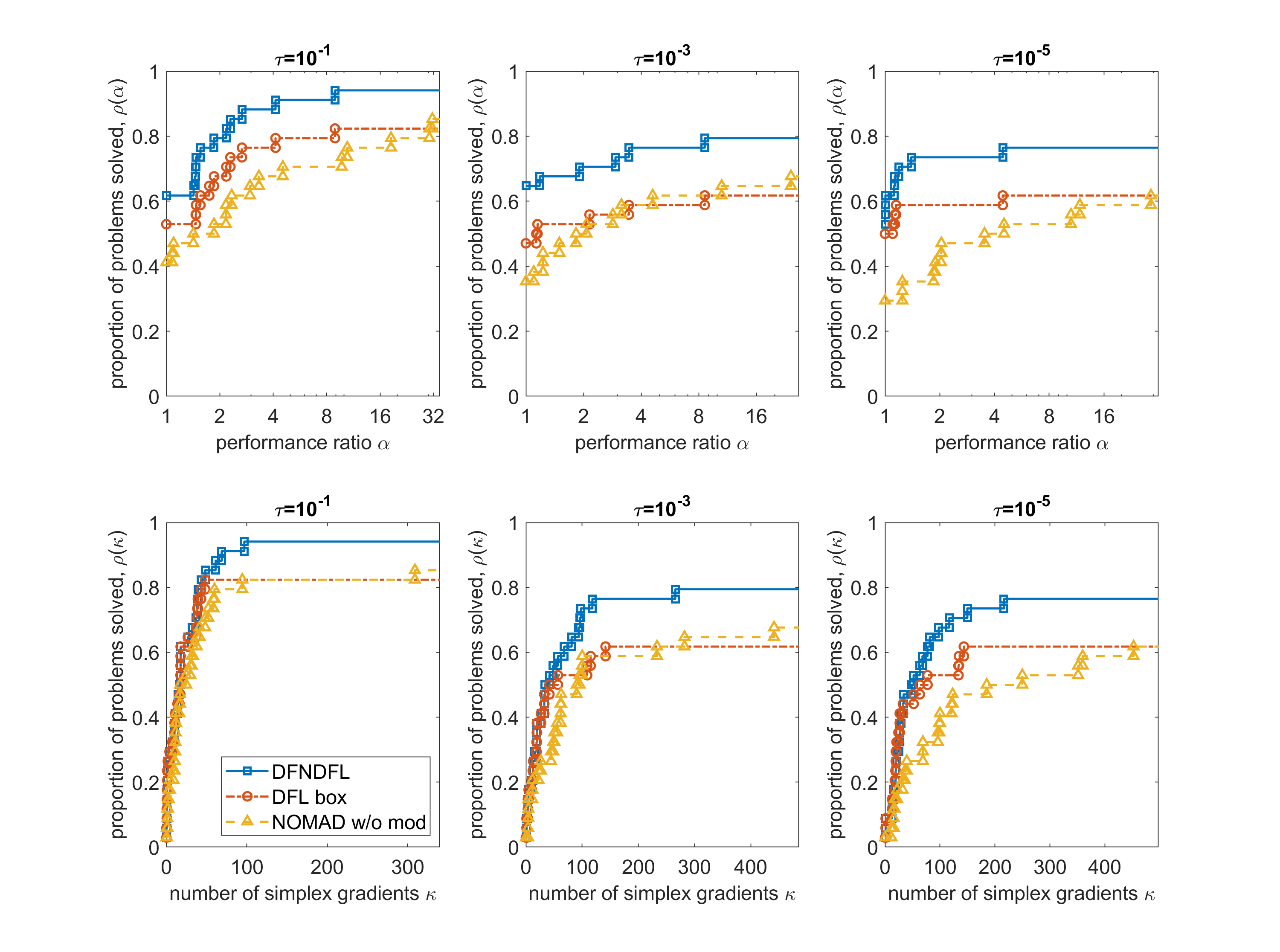}
 	\caption{Performance and data profiles for the comparison among DFNDFL, DFL {\em box}, NOMAD (not using models) on the 34 bound constrained problems. }\label{fig:box5000_1_e-1-3_perfprof_without_mod}\label{fig:box5000_1_e-1-3_dataprof_without_mod}
\end{figure}

Figure~\ref{fig:box5000_1_e-1-3_perfprof_without_mod} reports performance and data profiles related to the comparison of the algorithms that do not employ models, namely  DFNDFL, DFL {\em box} and NOMAD (without models). 
In this case, DFNDFL turns out to be the most efficient and reliable algorithm, regardless of the accuracy required. DFL {\em box} is more efficient and reliable than NOMAD when $\tau = 10^{-1}$ and $10^{-5}$. On the other hand, for $\tau = 10^{-3}$, we notice that DFL {\em box} is more efficient but less reliable than NOMAD (without models). It is worth noticing that the remarkable percentage of problems solved for $\alpha=1$ is an important result for DFNDFL since it shows that using more sophisticated directions than DFL {\em box} does not lead to a significant loss of efficiency. We also point out that the initial continuous and primitive search directions used by DFNDFL are the coordinate directions, which are the same as the ones employed in DFL {\em box}, thus leading to the same behavior of the algorithms in the first iterations. For each value of $\tau$, despite the remarkable efficiency, DFL {\em box} does not show a strong reliability, which is significantly improved by DFNDFL. 
\par

Next we compare DFNDFL against solvers that make use of sophisticated models to improve the search. In particular, we considered NOMAD (using models), RBFOpt and MISO. We point out that these three solvers exploit different kinds of models. In particular, NOMAD takes advantage of quadratic models whereas RBFOpt and MISO make use of radial basis function models. 

\begin{figure}[htbp]
	\centering
	\includegraphics[width=0.80\textwidth]{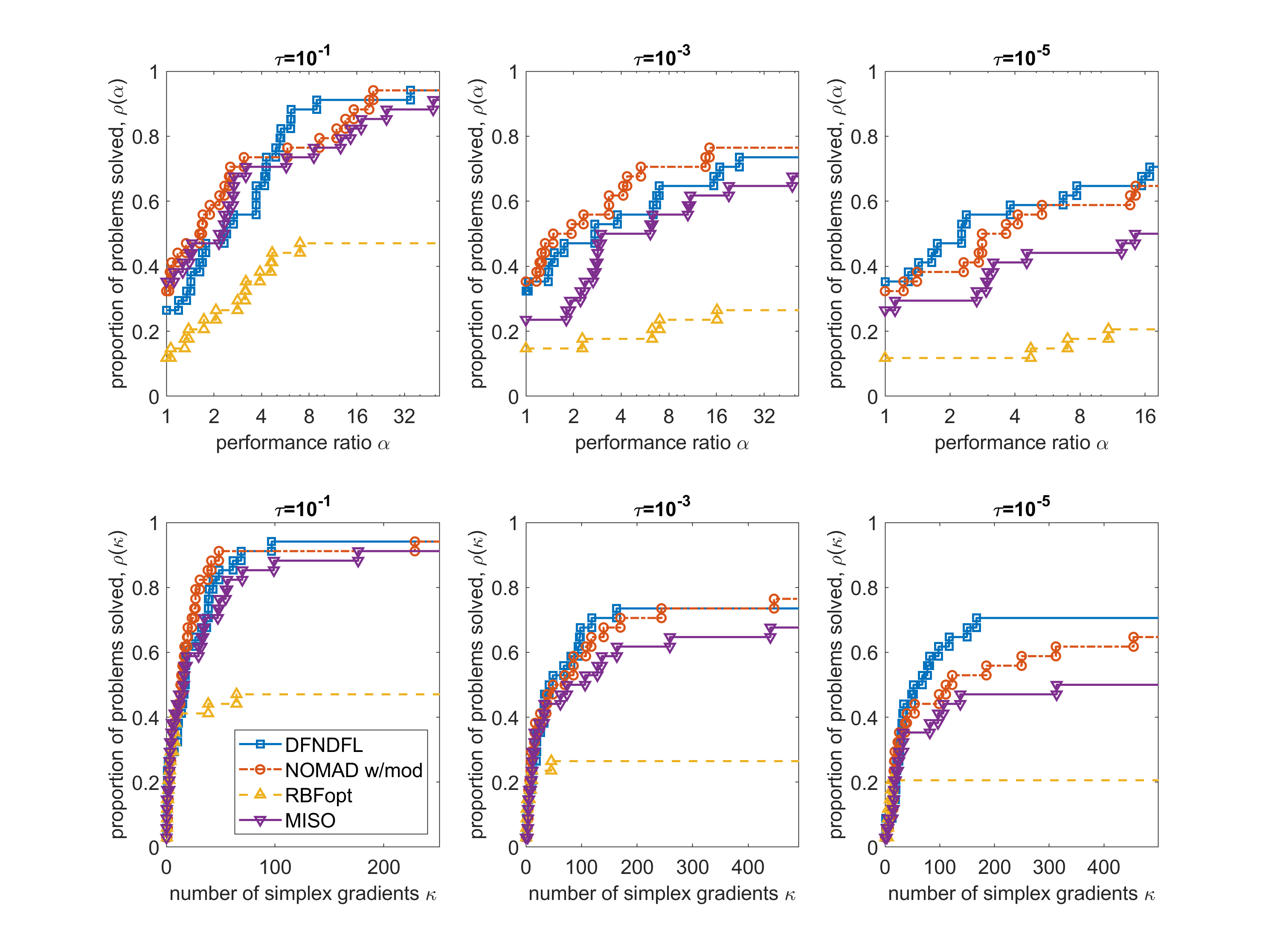}
	\caption{Performance and data profiles for the comparison among DFNDFL, NOMAD (using models), RBFOpt and MISO on the 34 bound constrained problems. }\label{fig:box5000_1_e-1-3_perfprof_with_mod}\label{fig:box5000_1_e-1-3_dataprof_with_mod}
\end{figure}

Figure~\ref{fig:box5000_1_e-1-3_perfprof_with_mod} reports performance and data profiles for the three considered levels of accuracy when DFNDFL is compared with the algorithms that make use of models.

From the performance and data profiles, it clearly emerges that RBFOpt is not very competitive with the other algorithms. Next, we can note that, despite its very relevant computational time required, MISO is outperformed by DFNDFL and NOMAD (w/ models) especially for medium ($\tau=10^{-3}$) and high ($\tau=10^{-5}$) accuracy levels. As concerns the comparison between DFNDFL and NOMAD (using models) the situation is somewhat more balanced even though for high accuracy DFNDFL is both more efficient and reliable than NOMAD.

\par
These numerical results highlight that DFNDFL has a remarkable efficiency and compares favorably to the state-of-the-art-solvers in terms of reliability, thus confirming and strengthening the properties of DFL {\em box} and providing a noticeable contribution to the derivative-free optimization solvers for bound constrained problems.
  
\subsubsection*{The nonlinearly constrained case}

The algorithms adopted for comparison in this case are DFL {\em gen} and the two versions of NOMAD (w/ and w/o models). We point out that the handling of the constraints in NOMAD is performed by the progressive/extreme barrier approach (by specifying the {\tt PEB} constraints type). We would like to highlight that we only used NOMAD and the constrained version of DFL {\em box}  (namely DFL {\em gen})  in this further comparison, due to their better performances in the bound constrained case and the explicit handling of  nonlinearly constrained problems. 

\begin{figure}[htbp]
	\centering
	\includegraphics[width=0.80\textwidth]{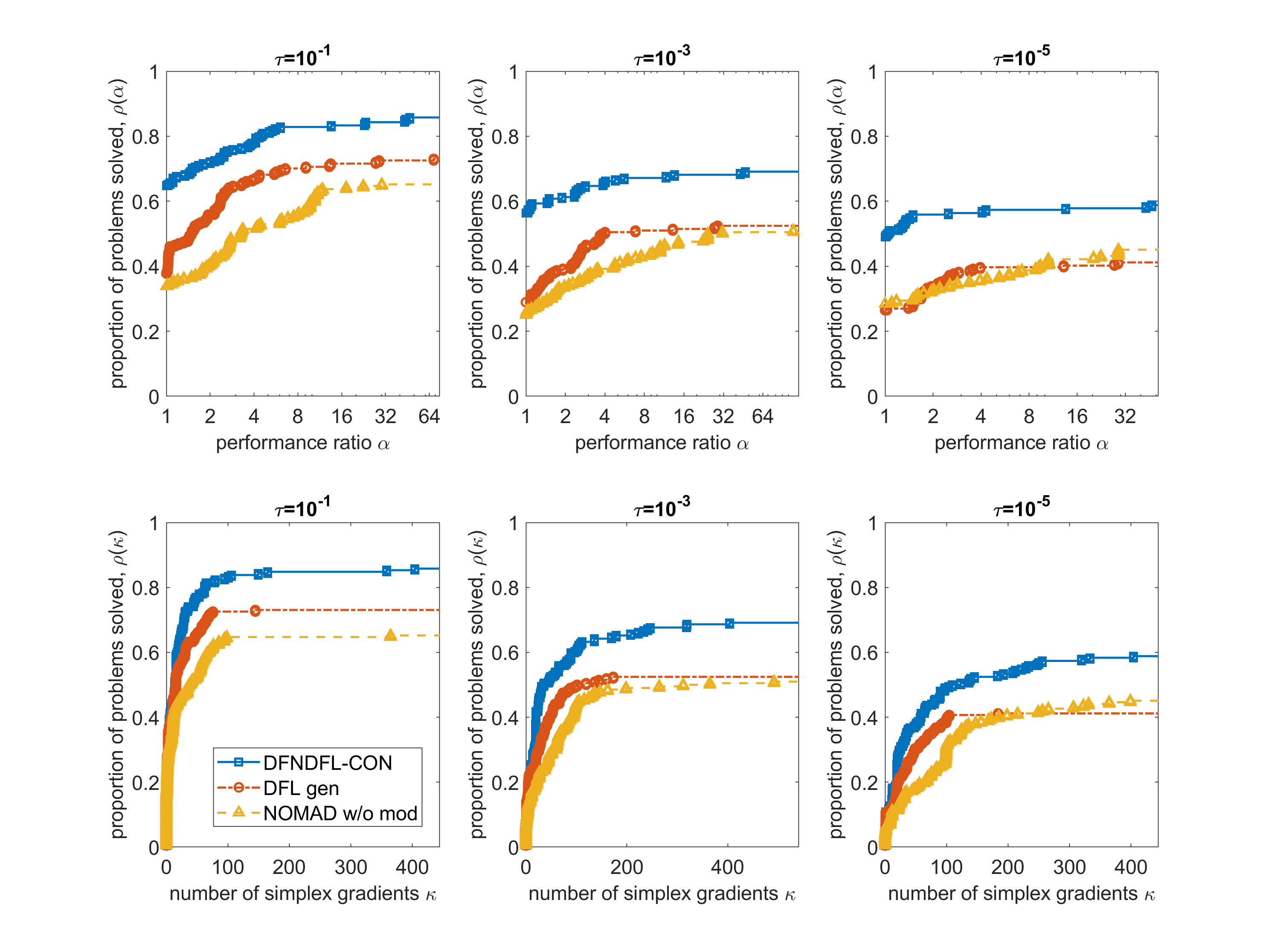} 
	\caption{Performance and data profiles for the comparison among DFNDFL--CON, DFL gen, NOMAD (w/o models) on the 204 nonlinearly constrained problems. }\label{fig:con_1_e-1-3_perfprof}\label{fig:con_1_e-1-3_dataprof}
\end{figure}
\begin{figure}[htbp]
	\centering
	\includegraphics[width=0.80\textwidth]{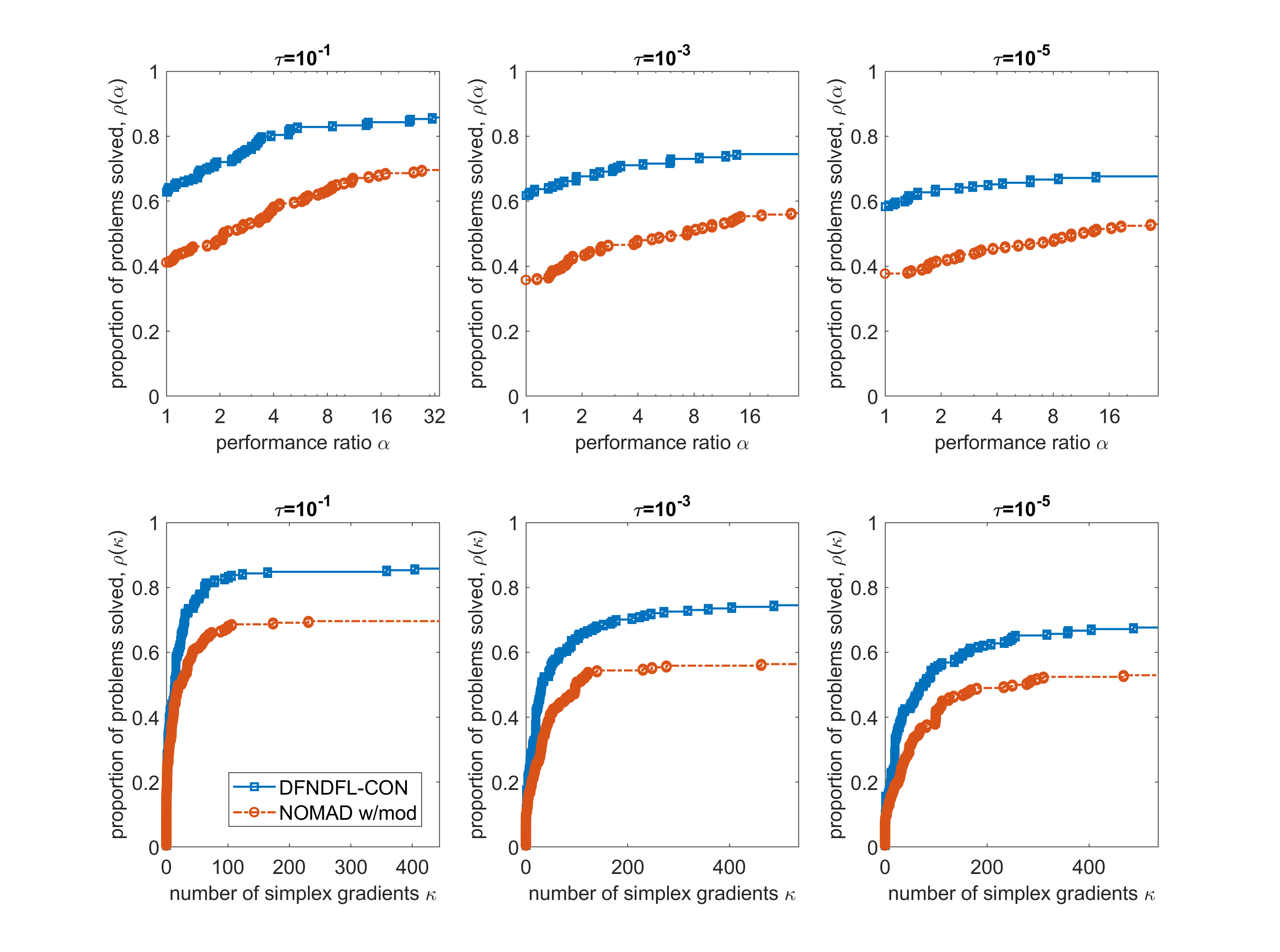} 
	\caption{Performance and data profiles for the comparison between DFNDFL--CON, and NOMAD (w/ models) on the 204 nonlinearly constrained problems. }\label{fig:dfndfl_nomadmod}
\end{figure}

Figure~\ref{fig:con_1_e-1-3_perfprof} reports  performance and data profiles, for the comparison of DFNDFL--CON, DFL {\em gen} and NOMAD (not using models). The figure quite clearly shows that DFNDFL--CON is the most efficient and reliable algorithm, and the difference with the other algorithms significantly grows as the level of accuracy increases. It is important to highlight that using the primitive directions allows our algorithm to improve the strategy of DFL, which only uses the set of coordinate directions. This results in a larger percentage of problems solved by DFNDFL--CON, even when compared with NOMAD (w/ mod). 



Finally, Figure \ref{fig:dfndfl_nomadmod} reports the comparison between DFNDFL--CON and NOMAD using models on the set of 204 constrained problems. Also in this case, it emerges that DFNDFL--CON is competitive with NOMAD both in terms of data and performance profiles. 

To conclude, these numerical results show that DFNDFL--CON has remarkable efficiency and reliability when compared to state-of-the-art-solvers.

\section{Conclusions}\label{sec:conclusion}
In this paper, new linesearch-based methods for mixed-integer nonsmooth optimization problems have been developed assuming that first-order information on the problem functions is not available. First, a general framework for bound constrained problems has been described. Then, an exact penalty approach has been  proposed and embedded into the framework for the bound constrained case, thus allowing to tackle the presence of nonlinear (possibly nonsmooth) constraints. Two different sets of directions are adopted to deal with the mixed variables. On the one hand, a dense sequence of continuous search directions is required to detect descent directions (those directions can form an arbitrarily narrow cone in the nonsmooth case). On the other hand, primitive  directions are employed to suitably explore the integer lattice (thus avoiding to get trapped into bad points).   

Numerical experiments have been performed on both bound and nonlinearly constrained problems. The results highlight that the proposed algorithms have good performances when compared with some state-of-the-art-solvers, thus providing a good tool for  handling the considered class of derivative-free optimization problems. 


\appendix 

\section{Technical results}

We first prove the equivalence between local minimum points and global minimum points of the two problems. Then, we prove that any Clarke stationary point of Problem \eqref{prob_Z} (according to Definition \ref{def:clarke_stationary_point}) is a stationary point for Problem \eqref{probconstr} (according to Definition \ref{def:stationary_point_constr}).


\begin{proposition}\label{th:no_stationary_points_out_feasible_region}	
	Let Assumption \ref{assmfcq} hold. A threshold value $\varepsilon^\star > 0$ exists such that the function $P(x;\varepsilon)$ has no Clarke stationary points in $(X\cap\mathcal{Z})\setminus\cal F$ for any $\varepsilon \in (0,\varepsilon^\star]$. 
\end{proposition}
\proof \quad We proceed by contradiction and assume that for any integer $k$, an $\varepsilon_k \leq 1/k$ and a stationary point for Problem \eqref{prob_Z} $x_k\in (X\cap\mathcal{Z})\setminus\cal F$ exist. Then, let us consider a limit point $\bar x\in (X\cap\mathcal{Z})\setminus\stackrel{\circ}{\cal F}$ of the sequence $\{x_k\}$ and, without loss of generality, let us call the corresponding subsequence as $\{x_k\}$ as well. Then, since $x_k\to\bar x$, the discrete variables remain fixed, i.e. $(x_k)_i = \bar x_i$, for all $i\in I^z$ and $k$ sufficiently large. Now, the proof continues by separately assuming that point $(i)$ or $(ii)$ of Assumption \ref{assmfcq} holds. 
\par
First we assume that point $(i)$ of Assumption \ref{assmfcq} holds at $\bar x$.
Therefore, a direction $\bar s\in D^c(\bar x)$ exists such that 
\[
\left(\xi^{g_i}\right)^\top \bar s <0 \ \mbox{ for all }\ \xi^{g_i}\in\partial_c g_i(\bar x) \mbox{ with } i\in\{ h \in \{1,\dots,m\}: g_h(\bar x) \geq 0\}.
\]
\commento{
}
In particular, it follows that
\begin{equation}\label{maxxi}
\max_{\tiny\begin{array}{c}\xi^{g_i}\in\partial_c g_i(\bar x)\\ i\in I(\bar x)\end{array}} \left(\xi^{g_i}\right)^\top \bar s = -\eta < 0,
\end{equation}
where $I(\bar x) = \left\{i\in\{1,\dots,m\}: g_i(\bar x) = \bar \phi(\bar x) \right\}$, $\bar \phi(x) = \max\left\{0,g_1(x),\dots,g_m(x)\right\}$ and $\eta$ is a positive scalar. Note that $\bar x\not\in\cal F$ implies $\bar \phi(\bar x) > 0$. 
\par
By \cite[Proposition 2.3]{lin.2009}, it follows that $\bar s\in D^c(x_k)$. Moreover, since $x_k$ satisfies the Definition \ref{def:stationary_point} of stationary point, we have 
that 
\begin{equation}\label{zetazero11}
0 \leq P^{Cl}(x_k;\varepsilon,\bar s) = \max_{\xi\in\partial_c P(x_k;\varepsilon)}\xi^\top \bar s.
\end{equation}
By \cite{clarke.1983}, we know that 
\begin{equation}
\label{eq:partial_1}
\partial_c P(x_k;\varepsilon) \subseteq \partial_c f(x_k) + \frac{1}{\varepsilon}\partial_c (\max\left\{0,g_1(x_k),\dots,g_m(x_k)\right\})
\end{equation}
and 
\begin{equation}
\label{eq:partial_2}
\partial (\max\left\{0,g_1(x_k),\dots,g_m(x_k)\right\}) \subseteq Co\left(\{\partial_c g_i(x_k): i\in I(x_k)\}\right),
\end{equation}
where $\beta^i\geq 0$ for all $i\in I(x_k)$ and $Co(A)$ denotes the convex hull of a set $A$ (see \cite[Theorem 3.3]{Rockafellar.1970}).
\par
Therefore, by \eqref{zetazero11}--\eqref{eq:partial_2}, $\xi^{f}_k \in \partial_c f(x_k)$, \, $\xi^{g_i}_k \in \partial_c g_i(x_k)$ and $\beta^i_k$ \mbox{ with } $i\in I(x_k)$ exist such that
\begin{equation}\label{derposxi}
\left(\xi^{f}_k + \frac{1}{\varepsilon_k}\sum_{i\in I(x_k)}\beta^i_k\xi^{g_i}_k\right)^\top \bar s \geq 0,
\end{equation}

$$
\sum_{i\in I(x_k)}\beta^i_k = 1 \ \text{ and } \ \beta^i_k \geq 0.
$$
Since $m$ is a finite number, there exists a subsequence of $\{x_k\}$ such that $I(x_k) = \bar I$. 

Then, recalling that $(x_k)_i=\bar x_i$, for all $i\in I^z$ and for $k$ sufficiently large, and since a locally Lipschitz continuous function has a generalized gradient which is locally bounded, it results that the sequences $\{\xi^{f}_k\}$ and $\{\xi^{g_i}_k\}$, with $i\in \bar I$, are bounded. Hence, we get that
\begin{subequations}\label{limconds}
	\begin{eqnarray}
	&& \label{limconds1}\xi^{f}_k\to \bar\xi^{f},\\
	&& \xi^{g_i}_k\to \bar\xi^{g_i} \ \mbox{for all}\ i\in \bar I,\\
	&& \beta^i_k\to\bar\beta^i \ \ \mbox{for all}\ i\in \bar I.
	\end{eqnarray}
\end{subequations}
Now the upper semicontinuity of $\partial_c f$ and $\partial_c g_i$, with $i\in\bar I$, at $\bar x$ (see Proposition 2.1.5
in \cite{clarke.1983}) implies that $\bar\xi^{f}\in \partial_c f(\bar x)$ and $\bar\xi^{g_i}\in \partial_c g_i(\bar x)$ for all
$i\in\bar I$.

The continuity of the problem functions guarantees that for $k$ sufficiently large
\[
\{i:g_i(\bar x)-\phi(\bar x) < 0\} \subseteq \{i: g_i(x_k) -\phi(x_k) < 0\}, 
\]
and, in turn, this implies that for $k$ sufficiently large
\[
\{i: g_i(x_k) -\phi(x_k) = 0\} = I(x_k) \subseteq I(\bar x) = \{i:g_i(\bar x)-\phi(\bar x) = 0\}. 
\]
Since $I(x_k) \subseteq I(\bar x)$, we have that
\begin{equation}\label{inclusione}
\bar I \subseteq I(\bar x). 
\end{equation}
Finally, for $k$ sufficiently large, \eqref{maxxi}, \eqref{limconds}, and \eqref{inclusione} imply 
\begin{equation}\label{derposxi2}
\left(\xi_k^{g_i}\right)^\top\bar s \leq -\frac{\eta}{2} \ \text{ for all } i\in\bar I,
\end{equation}
and (\ref{derposxi}), multiplied by $\varepsilon_k$, implies
\begin{equation} \label{eq:eq_3}
\left(\varepsilon_k\xi^{f}_k +\sum_{i\in\bar I}\beta^i_k\xi^{g_i}_k\right)^\top \bar s \geq 0.
\end{equation}
Equations \eqref{derposxi2} and \eqref{eq:eq_3} yield
\[
0 \leq \left(\varepsilon_k\xi^{f}_k +\sum_{i\in\bar I}\beta^i_k\xi^{g_i}_k\right)^\top \bar s \leq 
\left(\varepsilon_k\xi^{f}_k \right)^\top \bar s -\frac{\eta}{2}.
\]
which, by using \eqref{limconds}, gives rise to a contradiction when $\varepsilon_k \to 0$.

Now we assume that point $(ii)$ of Assumption \ref{assmfcq} holds at $\bar x$. Let  $\bar d \in D^z(\bar x)$ be the direction such that
\begin{equation}
	\label{eq:eq_dir}
	\sum_{i=1}^m \max\{0, g_i(\bar x+\bar d)\} < \sum_{i=1}^m \max\{0, g_i(\bar x)\}.
\end{equation} 
recalling that $(x_k)_i=\bar x_i$, for all $i\in I^z$ and for $k$ sufficiently large, we have that for $k$ sufficiently large $D^z(\bar x) = D^z(x_k)$, so that $\bar d\in D^z(x_k)$. 
By definition of stationary point and discrete neighborhood, we have
$$P(x_k; \varepsilon_k)\le P(x_k + \bar d; \varepsilon_k),\quad\quad \mbox{where}\ x_k + \bar d \in {\cal B}^z(x_k).$$
Hence,
$$f(x_k) + \frac{1}{\varepsilon_k} \sum_{i=1}^m \max\{0, g_i(x_k)\} \le f(x_k + \bar d) + \frac{1}{\varepsilon_k} \sum_{i=1}^m \max\{0, g_i(x_k+\bar d)\}.$$
Multiplying by $\varepsilon$ and considering that $\varepsilon_k \to 0$, if we take the limit for $k \to \infty$ we have that
$$\sum_{i=1}^m \max\{0, g_i(\bar x)\} \le \sum_{i=1}^m \max\{0, g_i(\bar x+\bar d)\}.$$
The latter equation is in contradiction with \eqref{eq:eq_dir}.
$\hfill\Box$ \\
\par
Now, we can prove that there exists
a threshold value $\bar \varepsilon$ for the penalty parameter such that, for any  $\varepsilon \in(0, \bar \varepsilon)$, any local minimum of the penalized problem is also a local minimum of the original problem.
\begin{proposition}\label{local_exactness}
	Let Assumptions \ref{assmfcq} hold. Given Problem \eqref{probconstr} and considering Problem \eqref{prob_Z}, a threshold value $\bar \varepsilon>0$ exists such that for every $\varepsilon \in(0, \bar \varepsilon]$, any
	local minimum point $\bar x$ of Problem \eqref{prob_Z} is also a local minimum of Problem \eqref{probconstr}.
\end{proposition}
\proof \quad 
Let $\bar x$ be any local minimum point of $P(x;\eps)$ on $X\cap {\cal Z}$.
By Corollary \ref{cor:clarke_stationary_point}, we have that $\bar x$ is also a Clarke stationary point. 
\par
Now Proposition \ref{th:no_stationary_points_out_feasible_region} implies that a threshold value $\varepsilon^\star > 0$ exists such that $\bar x \in {\cal F} \cap {\cal Z} \cap X$ for any $\varepsilon \in (0,\varepsilon^\star]$. Therefore, $P(\bar x;\varepsilon) = f(\bar x)$ which implies that $\bar x$ is also a local minimum point for Problem \ref{probconstr}.
$\hfill\Box$
\par\medskip\noindent

\begin{proposition}\label{global_exactness}
	Let Assumption \ref{assmfcq} hold. Given Problem \eqref{probconstr} and considering
	Problem \eqref{prob_Z}, a threshold value $\bar \varepsilon>0$ exists such that for every $\epsilon \in(0, \bar \varepsilon)$, any
	global minimum point $\bar x$ of Problem \eqref{prob_Z} is also a global minimum point of Problem \eqref{probconstr} and viceversa.
\end{proposition}
\proof \quad We start by proving that any global minimum point of Problem \eqref{prob_Z} is also a global minimum point of Problem \eqref{probconstr}. Proceeding by contradiction, let us assume that for any integer $k$ a positive scalar $\varepsilon_k < 1/k$ and a point $x_k$ exist such that $x_k$ is a global minimum point of $P(x_k;\varepsilon_k)$ but it is not a global minimum point of $f(x)$. If we denote as $\hat{x}$ a global minimum point of $f(x)$, we have that
\begin{equation} \label{eq:global_exactness_1}
P(x_k;\varepsilon_k) \le P(\hat{x};\varepsilon_k) = f(\hat{x}).
\end{equation}
Since $x_k$ are global minimum points, by Corollary \ref{cor:clarke_stationary_point} they are also stationary points of $P(x;\varepsilon_k)$ according to Clarke definition for the continuous variables. By Proposition \ref{th:no_stationary_points_out_feasible_region} there exists a threshold value $\varepsilon^\star > 0$ such that $x_k \in {\cal F} \cap X \cap {\cal Z}$ for any $\varepsilon_k \in (0,\varepsilon^\star]$. Therefore, $P(x_k;\varepsilon_k) = f(x_k)$. By \eqref{eq:global_exactness_1}, it follows that $f(x_k) \le f(\hat x)$, contradicting the assumption that $x_k$ is not a global minimum point of $f(x)$.
\par
Now we prove that any global minimum point $\bar x$ of Problem \eqref{probconstr} is also a global minimum point of Problem \eqref{prob_Z} for any $\varepsilon \in (0,\bar \varepsilon$). Since $\bar x \in {\cal F} \cap {\cal Z} \cap X$, we have that $P(\bar x;\varepsilon) = f(\bar x)$. By the previous proof, a global minimizer $x_\varepsilon$ of $P(x;\varepsilon)$ is feasible for Problem \eqref{probconstr}, hence $P(x_\varepsilon;\varepsilon) = f(x_\varepsilon)$. Furthermore, it is also a global minimum point of Problem \eqref{probconstr}, thus we have $f(x_\varepsilon) = f(\bar x)$. Therefore, since $P(x_\varepsilon;\varepsilon) = f(\bar x)$, $\bar x$ is also a global minimum point of $P(x;\varepsilon)$. 
$\hfill\Box$
\par\medskip\noindent

In order to give stationarity results for Problem \eqref{prob_Z}, we have the following proposition.

\begin{proposition}\label{equiv_feasible}
	Let Assumption \ref{assmfcq} hold. For any $\varepsilon >0$, every stationary point $\bar x$ of Problem \eqref{prob_Z} according to Clarke, such that $\bar x\in{\cal F} \cap {\cal Z} \cap X$, is also a stationary
	point of Problem \eqref{probconstr}.
\end{proposition}
\proof \quad
Since $\bar x$ is, by assumption, a stationary point of Problem \eqref{prob_Z} according to Clarke (see Corollary \ref{cor:clarke_stationary_point}), then we have by definition of Clarke stationarity that
for all $s\in D^c(\bar x)$, 
\begin{align}\label{eq:stationarity_cont} P^{Cl}(\bar x;\varepsilon) = \max \left \{ \xi^\top s\, : \, \xi\in \partial_c P(\bar x;\varepsilon)\right \}\geq 0,
\end{align}
and 
\begin{align}\label{eq:stationarity_discr}
P(\bar x;\varepsilon)\le P(x;\varepsilon)\quad\mbox{for all}\ x\in {\cal B}^z(\bar x).
\end{align} 
Hence, by \ref{eq:stationarity_cont}, there exists $\xi_s\in\partial_c  P(\bar x;\varepsilon)$ such that $(\xi_s)^\top s \geq 0$  for all $s\in D^c(\bar x)$.
Now, we recall that 
$$
\partial_c P(x;\varepsilon) \subseteq \partial_c f(x) + \frac{1}{\varepsilon}\sum_{i\in I(x)} \beta_i\partial_c g_i(x),
$$
for some $\beta_i$, with $i\in I(x)$, such that $\sum_{i\in I(x)}\beta_i=1$ and $\beta_i\geq 0$ for all $i\in I(x)$. Hence,
we have that $\xi_s \in \partial_c f(\bar x) + \frac{1}{\varepsilon}\sum_{i\in I(\bar x)} \beta_i\partial_c g_i(\bar x)$. Then,
denoting $\lambda_i = \beta_i/\varepsilon$ with $i\in I(\bar x)$, and assuming $\lambda_i = 0$ for all $i \notin I(\bar x)$, we can write, for all $s \in D^c(\bar x)$, 
\begin{align} 
&\max \left \{\xi^\top s\, : \, \xi\in \partial_c f(\bar x) + \displaystyle\sum_{i=1}^{m} \lambda_i\partial_c g_i(\bar x) \right \}\geq 0, \label{eq:kkt_stationarity_1} \\
&(\lambda)^T g(\bar x) = 0 \text{ and } \lambda\geq 0. \label{eq:kkt_stationarity_2}
\end{align}
Recalling that $\bar x$ is feasible for Problem \eqref{probconstr}, by \eqref{eq:stationarity_discr} we have
\begin{equation}\label{eq:stationarity_discr_f}
f(\bar x)\le f(x)\quad\mbox{for all}\ x\in {\cal B}^z(\bar x),
\end{equation}
Considering that $\bar x\in{\cal F} \cap {\cal Z} \cap X$, (\ref{eq:kkt_stationarity_1}), (\ref{eq:kkt_stationarity_2}),  and (\ref{eq:stationarity_discr_f}) prove that $\bar x$ is a KKT stationary point for Problem \eqref{probconstr}, thus concluding the proof. 
$\hfill\Box$ \\
\par
\begin{proposition}\label{equivalenza}
	Let Assumption \ref{assmfcq} hold. Then, a threshold value $\varepsilon^\star >0$ exists such that, for every $\varepsilon\in (0,\varepsilon^\star]$, every stationary point $\bar x$ 
	of Problem \eqref{prob_Z} is stationary (according to Definition \ref{def:stationary_point_constr}) for Problem (\ref{probconstr}).
\end{proposition}
\proof \quad Since $\bar x$ is stationary for Problem \eqref{prob_Z}, we have 
by Definition \ref{def:stationary_point} that
\begin{equation}\label{stat_cond3}
P^\circ(\bar x; \varepsilon, s)\geq 0 \ \text{ for all } \ s\in D^c(\bar x),
\end{equation}
and
\begin{align}\label{eq:stationarity_discr_2}
P(\bar x;\varepsilon)\le P(x;\varepsilon)\quad\mbox{for all}\ x\in {\cal B}^z(\bar x).
\end{align}
Then, by Definitions  \ref{def:clarke_directional_derivative} and \ref{def:clarke-jahn_directional_derivative}, we have that
$$
\limsup_{\footnotesize\begin{array}{l}y_c\to x_c, y_z = x_z, t\downarrow 0\end{array}} \frac{P(y+t s;\varepsilon) - P(y;\varepsilon)}{t} = 
P^{Cl}(\bar x;\varepsilon,s)\geq P^\circ(\bar x; \varepsilon,s)
\ \text{for all } \ s \in D^c(\bar x).
$$
By \eqref{stat_cond3}, it follows that
$$
P^{Cl}(\bar x;\varepsilon,s)\geq 0
\ \text{ for all } \ s \in D^c(\bar x).
$$
The proof follows by considering Propositions \ref{th:no_stationary_points_out_feasible_region}, \ref{equiv_feasible} and (\ref{eq:stationarity_discr_2}). 
$\hfill\Box$


\end{document}